\ifx\formloaded\relax   \else \let\formloaded=\relax \fi

\hsize=16cm
\vsize=24cm
\normalbaselineskip=13pt\normallineskiplimit=.1pt\normallineskip.1pt
\normalbaselines
\setbox\strutbox=\hbox{\vrule height8.9pt depth4pt width0pt}
\parskip=\smallskipamount 
\multiply \smallskipamount by 4 \divide\smallskipamount by 3
\multiply \medskipamount by 4 \divide\medskipamount by 3
\multiply \bigskipamount by 4 \divide\bigskipamount by 3


\ifx\CWIloaded\undefined \let\CWIloaded\relax \else   \fi

\ifx\topglue\undefined \def\topglue{\nointerlineskip\vglue-\topskip\vskip}\fi

\ifx\formloaded\undefined 
  \vsize=24.5 cm \advance \vsize -24pt 
  \hsize=16 cm
  \let\formloaded\relax 
\fi


\catcode`_=11 

\def\let_next{\let\next=}
{\obeylines\gdef\default_address%
 {Universit\'e de Poitiers, D\'epartement de Math\'ematiques,
  SP2MI, T\'el\'eport~2, BP~179, \ 86960 Futuroscope Cedex, France}}

\newtoks\title \newtoks\MSCcodes \newtoks\CRcodes \newtoks\keywords
\newtoks\titlenote
\newif\if_started

 3

\font\seventeenssbf=cmssbx10 scaled \magstep 3

\font\twelvess=cmss12

\font\tenss=cmss10
\font\tenssit=cmssi10
\font\tenssbf=cmssbx10
\def\sans
{\def\rm{\fam0\tenss}\def\it{\fam\itfam\tenssit}\def\bf{\fam\bffam\tenssbf}\rm
}


\font\ninerm=cmr9 \font\nineit=cmti9 \font\ninebf=cmbx9 \font\ninesl=cmsl9
\font\ninett=cmtt9 \font\ninei=cmmi9 \font\ninesy=cmsy9
\font\niness=cmss9 \font\ninessit=cmssi9 \font\ninessbf=cmssbx10 at 9 pt
\skewchar\ninei='177 \skewchar\ninesy='60

\def\ninept
{\def\rm{\fam0\ninerm}\def\it{\fam\itfam\nineit}\def\bf{\fam\bffam\ninebf}%
 \def\sl{\fam\slfam\ninesl}\def\tt{\fam\ttfam\ninett}\rm
 \def\sans{\def\rm{\fam0\niness}\def\it{\fam\itfam\ninessit}%
           \def\sl{\fam\slfam\ninessit}\def\bf{\fam\bffam\ninessbf}\rm}%
 \textfont0\ninerm \textfont1\ninei \textfont2=\ninesy
 \textfont\itfam\nineit \textfont\slfam\ninesl
 \textfont\bffam\ninebf \textfont\ttfam\ninett
 \normalbaselineskip=11pt \normalbaselines
 \setbox\strutbox=\hbox{\vrule height 8 pt depth 3 pt width 0 pt}%
}


\long\def\centerpar#1{\begingroup
  \advance\leftskip by 0pt plus 1 fil minus 3pt \rightskip=\leftskip
	\parskip=0pt \parindent=0pt \parfillskip=0pt \linepenalty=100
  \noindent\ignorespaces#1\par\endgroup}

\def\good_break{\unskip\penalty-10\ \ignorespaces}

\catcode`@=11
\def\vfootnote#1%
{\insert\footins\bgroup \ninept
 \interlinepenalty=\interfootnotelinepenalty
 \splittopskip=\ht\strutbox 
 \splitmaxdepth=\dp\strutbox 
 \floatingpenalty=20000 
 \leftskip=1pc \rightskip=1pc \spaceskip=0pt \xspaceskip=0pt
 \parindent=0pt \parskip=3pt minus 1pt \parfillskip=0pt plus 1 fil
 \textindent{#1}\footstrut\futurelet\next\f@@t
}
\catcode`@=12

\let\_authors\empty 

%

\def\Author#1{\def\this_author{#1}%
  \ifx\this_author\_eq \let_next\Author 
  \else \let_next\let_next \afterassignment\_author 
  \fi \next}
\def\_author
 {\let\this_addr=\default_address
  \ifx\next\address \let_next\get_address
  \else \ifx\next\email \let_next\get_email
  \else \ifx\next\URL \let_next\get_URL
  \else \add_author
  \fi\fi\fi \next
}

\def\address#1{\errmessage
 {\string\address\space should only follow \string\Author}}
\def\email#1{\errmessage
 {\string\email\space should only follow \string\Author}}
\def\URL#1{\errmessage
 {\string\URL\space should only follow \string\Author}}

\def\get_address{\begingroup\obeylines\get_addr}
\def\get_addr#1%
{\endgroup\def\this_addr{#1}%
  \ifx\this_addr\_cr \let_next\get_address 
  \else\ifx\this_addr\_eq \let_next\get_address 
       \else \let_next\let_next \afterassignment\get_ad
  \fi\fi \next
}
\def\get_ad
{\ifx\next\email \let_next\get_email
 \else \ifx\next\URL \let_next\get_URL
 \else \add_author
 \fi\fi \next
}

\def\get_email#1%
{\def\this_email{#1}%
  \ifx\this_email\_eq \let_next\get_email
  \else \edef\this_addr{\this_addr\endgraf{\tt #1}}%
    \let_next\let_next \afterassignment\after_em
  \fi \next
}
\def\after_em{\ifx\next\URL \let_next\get_URL \else \add_author \fi \next}
\def\get_URL#1%
{\def\this_URL{#1}%
 \ifx\this_URL\_eq \let_next\get_URL
 \else 
  {\def~{\char "7E }%
   \xdef\this_addr{\this_addr\endgraf{\noexpand\tt URL: http://#1}}%
  }\let_next\add_author
 \fi \next
}

\def\_eq{=}{\obeylines\gdef\_cr{^^M}}

\def\add_author 
{\toks0=\expandafter{\_authors}%
 \xdef\_authors
 {\the\toks0 \noexpand\\{\this_author\noexpand\_addr{\this_addr}}}%
}


\def\start_art
{\vfil\eject \ifodd\pageno\else \shipout\vbox to \vsize{}\advancepageno\fi

 \topglue 10mm plus 3 mm minus 2mm
 \centerpar{\let\\=\good_break\seventeenssbf \baselineskip=22pt \the\title}%
 \ifx \_authors\empty \message{Warning: no authors specified.}\else
   \bigskip\centerpar
    {\twelvess
     \def\par{\ifhmode\endgraf\else\medskip\fi}\obeylines
     \def\_addr##1{\smallskip{\ninessit\let\\=\par\noindent##1\endgraf}}%
     \def\\##1{\bigskip##1}\_authors
    }\bigskip
 \fi
 \_startedtrue
}


\outer\def\abstract
{\if_started \errmessage{Abstract must be placed at start of article.}%
 \else \start_art \bigskip \begingroup \ninept \sans \narrower
   \leftline{\hskip\leftskip ABSTRACT}\nobreak\smallskip
	 \everypar{\everypar{}\setbox0=\lastbox}%
 \fi
}

\outer\def\maintext
{\if_started\medskip\endgroup\else\start_art\fi
 \if &\the\MSCcodes\the\CRcodes\the\keywords\the\titlenote&\else
  {\narrower\ninept\sans
   \if &\the\MSCcodes&\else
     \hangindent=4em \noindent
     {\it 1991 Mathematics Subject Classification:\/} \the\MSCcodes.
     \par\fi
   \if &\the\CRcodes&\else
     \hangindent=4em \noindent
     {\it 1991 Computing Reviews Subject Classification:\/} \the\CRcodes.
     \par\fi
   \if &\the\keywords&\else
     \hangindent=2em \noindent{\it Keywords and Phrases:\/} \the\keywords.
     \par\fi
   \if &\the\titlenote&\else
     {\def\\{\unskip\hfil\break\ignorespaces}%
      \hangindent=2em \noindent{\it Note:\/} \the\titlenote.
      \par}\fi
  }\bigskip
 \fi
}

\def\sticksection 
{\vbox\bgroup\parskip=0pt \everypar{\egroup\noindent}}

\def\newreport
{\vfil \eject \secno=0\subsecno=0\proclno=0 \eqnumber=0 \pageno=1 \mark{}
 \title={}\MSCcodes={}\CRcodes={}\keywords={}\titlenote={}
 \let\_authors\empty
 \_startedfalse
}

\catcode`_=8 
\input marxmax

\refindent=3.2pc



\newif\ifnarrow
\ifdim \hsize<14.5cm \narrowtrue \else\narrowfalse \fi

\def\maybesmall{\ifnarrow\smallsquares\fi}

\def\ol#1
{\vbox{\vskip0.3pt \hrule \vskip 0.8pt \hbox{$#1$}}}

\def\ul#1
{\vtop{\hbox{$#1$}\vskip 0.8pt \hrule \vskip0.3pt }}

\font\frak=eufm10
\font\csc=cmcsc10
\font\URLfont=cmtex10 
\hyphenchar\URLfont=-1 

\mathorddef\Part=P2 
\mathorddef\A=A0 \mathorddef\B=B0 
\def\GL{\boldmath{GL}}
\def\Ylat{(\Part,{\subset})}
\def\upto#1{[#1]}
\opdef\Pic
\altopdef\Tab{SST}
\opdef\ST
\opdef\LR
\opdef\JdT
\opdef\wt
\mathorddef\T=T2 
\mathorddef\Sk=S2 
\def\SFn{\Lambda_n}
\def\Can#1{\boldmath1_{#1}}
\def\Ylm{Y(\\/\mu)}
\def\bT{\bar T}
\def\LRt#1,#2{\LR(#1,#2)}
\def\JdTt#1,#2{\Tab(#1)^{\slid{#2}}}
\mathorddef\Rob=R2 
\def\StT#1{\ST(#1)}
\def\snu{^{\slid\nu}}
\def\Uqgln{U_q(\hbox{\frak gl}_n)}
\def\leqr{\leq_r}
\def\leqKan{\leq_{\romath K}}

\def\wKan{w_{\romath K}}
\def\wSem{w_{\romath S}}

\mathdef\slid=Rel \triangleright
\mathdef\dils=Rel \triangleleft

\long\def\quotation#1
{{\smallskip\narrower\it\noindent\llap{``\/}#1\unskip''\smallskip}\noindent}

{\catcode`\.=\active\catcode`\/=\active\catcode`\_=\active
 \gdef\URL
 {\begingroup\def~{\char`\~}%
  \catcode`\.=\active\catcode`\/=\active\catcode`\_=\active
  \let_=\_%
  \def.{\char`\.\penalty\hyphenpenalty}\def/{\char`\/\penalty\hyphenpenalty}%
  \URLfont\URLnext
}}
\def\URLnext#1{http://#1\endgroup}

\hyphenation
{semi-standard Yamano-uchi mono-mial mono-mials inter-related combi-na-torics
 combi-na-torial algo-rithm valid-ity}

\ifx\Tableauxloaded\relax   \else \let\Tableauxloaded\relax \fi


\newcount\cols
{\catcode`,=\active\catcode`|=\active
 \gdef\Young(#1){\hbox{$\vcenter
 {\mathcode`,="8000\mathcode`|="8000
  \def,{\global\advance\cols by 1 &}%
  \def|{\cr
        \multispan{\the\cols}\hrulefill\cr
        &\global\cols=2 }%
  \offinterlineskip\everycr{}\tabskip=0pt
  \dimen0=\ht\strutbox \advance\dimen0 by \dp\strutbox
  \halign
   {\vrule height \ht\strutbox depth \dp\strutbox##
    &&\hbox to \dimen0{\hss$##$\hss}\vrule\cr
    \noalign{\hrule}&\global\cols=2 #1\crcr
    \multispan{\the\cols}\hrulefill\cr%
   }
 }$}}
 \gdef\Skew(#1:#2){\hbox{$\vcenter
 {\mathcode`,="8000\mathcode`|="8000
  \dimen0=\ht\strutbox \advance\dimen0 by \dp\strutbox
  \def\boxbeg{\vbox
    \bgroup\hrule\kern-0.4pt\hbox to\dimen0\bgroup\strut\vrule\hss$}%
  \def\boxend{$\hss\egroup\hrule\egroup}%
  \def,{\boxend\boxbeg}%
  \def|##1:{\boxend\vrule\egroup\nointerlineskip\kern-0.4pt
    \moveright##1\dimen0\hbox\bgroup\boxbeg}%
  \def\\##1\\##2:{\boxend\vrule\egroup\nointerlineskip\kern-0.4pt
    \kern ##1\dimen0\moveright##2\dimen0\hbox\bgroup\boxbeg}%
  \moveright#1\dimen0\hbox\bgroup\boxbeg#2\boxend\vrule\egroup
 }$}}
}


\font\sevenit=cmti7 \font\fiveit=cmti5
\scriptfont\itfam=\sevenit \scriptscriptfont\itfam=\fiveit

\def\smallsquares
{\textfont0=\scriptfont0 \scriptfont0=\scriptscriptfont0
 \textfont1=\scriptfont1 \scriptfont1=\scriptscriptfont1
 \textfont\itfam=\scriptfont\itfam \scriptfont\itfam=\scriptscriptfont\itfam
 \textfont\bffam=\scriptfont\bffam \scriptfont\bffam=\scriptscriptfont\bffam
 \setbox0=\hbox{$($}
 \setbox\strutbox=\hbox{\vrule width 0pt height\ht0 depth\dp0 }
}

\def\bigsquares
{\setbox\strutbox=\hbox
 {\vrule width 0pt height1.3\ht\strutbox depth1.3\dp\strutbox}
}

\newbox\hpair \newbox\vpair
\def\installboxes
{\dimen0=\ht\strutbox \advance\dimen0 by \dp\strutbox
 \setbox0=\vbox{\hrule width \dimen0}
 \setbox1=\rlap{\strut\vrule}
 \setbox\hpair=\rlap{\raise\ht\strutbox\copy0 \lower\dp\strutbox\copy0}
 \setbox\vpair=\rlap{\raise\dimen0\copy1 \kern\dimen0 \copy1}
}

\newcount\ribbonlength
\newcount\joff \newcount\jmax \newcount\imid \newcount\jmid
\def\ribbon #1,#2:#3;#4
{\setbox2=\hbox
 {\count0=0 \count1=0
  \copy1 \rlap{\lower\dp\strutbox\copy0}\if|#4|\else\ribs#4 \fi
  \raise\count1\dimen0\rlap
    {\kern\count0\dimen0 \raise\ht\strutbox\copy0 \copy1\raise0.4pt\copy1}%
  \global\joff=\count0
  \divide\dimen0 by 2 \raise\imid\dimen0\rlap{\kern\jmid\dimen0
    \hbox to 2\dimen0{\hfil$#3$\hfil}}%
 }%
 \lower#1\dimen0\rlap{\kern#2\dimen0\box2}%
 \advance\joff by #2\relax \ifnum \joff>\jmax \global\jmax=\joff \fi
}
\def\ribs#1#2
{\raise\count1\dimen0\rlap
   {\kern\count0\dimen0 \copy\ifnum #1=0 \hpair \else \vpair \fi}%
 \advance\ribbonlength by -2
 \ifnum \ribbonlength=0 \jmid=\count0 \imid=\count1 \fi
 \advance\count#1 by 1
 \ifnum \ribbonlength=1 \jmid=\count0 \imid=\count1 \fi
 \ifnum \ribbonlength>-1 \ifnum \ribbonlength<2
   \advance \jmid by \count0 \advance \imid by \count1
 \fi\fi
 \if |#2|\else \ribs#2 \fi
}

\def\rtab#1
{\installboxes \ribbonlength=#1 \global\jmax=0
 \setbox0=\hbox{\aftergroup\endrtab \aftergroup}}
\def\endrtab{\advance\dimen0 by \jmax\dimen0 \wd0=\dimen0 \vcenter{\box0}}

\def\arrow#1,#2:#3
{\lower#1\dimen0\rlap{$\kern#2.5\dimen0
 \if r#3 \,\rightarrow \else
 \if l#3 \llap{$\leftarrow\,$} \else
 \if u#3 \setbox0=\hbox to 0pt{\hss$\uparrow$\hss}\ht0=0pt 
         \raise.5\dimen0 \box0 \else
 \if d#3 \setbox0=\hbox to 0pt{\hss$\downarrow$\hss}\dp0=-2pt
         \lower.5\dimen0 \box0 \else
 \if i#3 \setbox0=\hbox to 0pt{\hss$\nwarrow\,$}\ht0=0pt
         \raise.5\dimen0 \box0 \else
 \if o#3 \setbox0=\hbox to 0pt{$\,\searrow$\hss}\dp0=-2pt
         \lower.5\dimen0 \box0 \else
 \fi\fi\fi\fi\fi\fi
 $}%
}

\readrefs

\Author={Marc A. A. van Leeuwen}
\email={maavl@mathlabo.univ-poitiers.fr}
\URL{wwwmathlabo.univ-poitiers.fr/~maavl/}

\title={The Littlewood-Richardson rule, \\
        and related combinatorics}
\MSCcodes={05E05, 05E10}


\abstract
An introduction is given to the Littlewood-Richardson rule, and various
combinatorial constructions related to it. We present a proof based on tableau
switching, dual equivalence, and coplactic operations. We conclude with a
section relating these fairly modern techniques to earlier work on the
Littlewood-Richardson rule.

\maintext

\let\endofchapter=\par

\secno=-1


\newsection Introduction.

The Littlewood-Richardson rule is a combinatorial rule describing the
multiplication of Schur polynomials; it was first formulated
in~\ref{Littlewood Richardson}, but its general validity remained unproved for
several decades. The various proofs that have been given since have created a
rich combinatorial theory, with many interrelated constructions, including the
Robinson-Schensted correspondence, and \foreign{jeu de taquin}. We describe
several such constructions, and use them to prove the rule; we are not however
narrowly focussed on this proof, and discuss several topics that are not used
in it. We make use of certain ``modern'' (post-1980) notions, while we do not
treat some other notions that figure prominently in many other proofs of the
Littlewood-Richardson rule, and are well documented elsewhere (e.g.,
\ref{Fulton}, \ref{Fomin appendix}); specifically, we focus on tableau
switching (which includes jeu de taquin), dual equivalence, and coplactic
operations
\unskip\footnote*
{Arguably coplactic operations are not at all new: even if not considered in
 isolation, they do occur in various forms and contexts, as a part of larger
 constructions; see \ref{Littlewood Richardson}, \ref{Robinson}, \ref{de
 Bruijn Tengbergen Kruyswijk}, \ref{Greene Kleitman}, \ref{Thomas Robinson},
 \ref{Lascoux Schutzenberger}.
},
but do not introduce Schensted's algorithm, Knuth equivalence, or Greene's
poset invariant. We do not really define or use Zelevinsky's pictures, but
they are mentioned in several places, and have inspired much of our work. At
the end, we give an overview of earlier work on the subject, to help clarify
its relation to our approach.

The proof we give does not require many technical verifications; moreover, we
believe that our main theorem~3.3.1 makes the correspondences considered more
transparent. We do prove in detail some properties of the constructions that
are so basic that they might have been left to the reader; this is because we
feel that it is often these ``low-level'' properties that really explain why
the more significant theorems work. Most facts presented in this paper are
known (at least to experts), although for some it is hard to find a published
reference. Nonetheless the global structure of our proof, and theorem~3.3.1,
seem to be new, even if the latter is related to the known fact that
``coplactic operations are compatible with plactic equivalence'' (which in
fact motivates their name).

The remainder of this paper is organised in four sections: in \Sec1 we
formulate the Littlewood-Richardson rule; in \Sec2 we define tableau switching
and derive an expression for Littlewood-Richardson coefficients in terms of
jeu de taquin; in \Sec3 we define coplactic operations and establish our main
theorem, which implies the Littlewood-Richardson rule; finally in \Sec4 we
comment on earlier work.

A word on notation: we always start indexing at~$0$; in particular this
applies to parts of partitions, and rows, columns, and entries of tableaux. We
define $\upto n=\setof i\in\N:i<n\endset$ for $n\in\N$.

 \par

\ifx\macrosloaded\undefined
  \input marxmax
  

\newif\ifnarrow
\ifdim \hsize<14.5cm \narrowtrue \else\narrowfalse \fi

\def\maybesmall{\ifnarrow\smallsquares\fi}

\def\ol#1
{\vbox{\vskip0.3pt \hrule \vskip 0.8pt \hbox{$#1$}}}

\def\ul#1
{\vtop{\hbox{$#1$}\vskip 0.8pt \hrule \vskip0.3pt }}

\font\frak=eufm10
\font\csc=cmcsc10
\font\URLfont=cmtex10 
\hyphenchar\URLfont=-1 

\mathorddef\Part=P2 
\mathorddef\A=A0 \mathorddef\B=B0 
\def\GL{\boldmath{GL}}
\def\Ylat{(\Part,{\subset})}
\def\upto#1{[#1]}
\opdef\Pic
\altopdef\Tab{SST}
\opdef\ST
\opdef\LR
\opdef\JdT
\opdef\wt
\mathorddef\T=T2 
\mathorddef\Sk=S2 
\def\SFn{\Lambda_n}
\def\Can#1{\boldmath1_{#1}}
\def\Ylm{Y(\\/\mu)}
\def\bT{\bar T}
\def\LRt#1,#2{\LR(#1,#2)}
\def\JdTt#1,#2{\Tab(#1)^{\slid{#2}}}
\mathorddef\Rob=R2 
\def\StT#1{\ST(#1)}
\def\snu{^{\slid\nu}}
\def\Uqgln{U_q(\hbox{\frak gl}_n)}
\def\leqr{\leq_r}
\def\leqKan{\leq_{\romath K}}

\def\wKan{w_{\romath K}}
\def\wSem{w_{\romath S}}

\mathdef\slid=Rel \triangleright
\mathdef\dils=Rel \triangleleft

\long\def\quotation#1
{{\smallskip\narrower\it\noindent\llap{``\/}#1\unskip''\smallskip}\noindent}

{\catcode`\.=\active\catcode`\/=\active\catcode`\_=\active
 \gdef\URL
 {\begingroup\def~{\char`\~}%
  \catcode`\.=\active\catcode`\/=\active\catcode`\_=\active
  \let_=\_%
  \def.{\char`\.\penalty\hyphenpenalty}\def/{\char`\/\penalty\hyphenpenalty}%
  \URLfont\URLnext
}}
\def\URLnext#1{http://#1\endgroup}

\hyphenation
{semi-standard Yamano-uchi mono-mial mono-mials inter-related combi-na-torics
 combi-na-torial algo-rithm valid-ity}

  \writelab{}
  \def\jobname{lrr}
\readrefs
  \let\endofchapter=\bye
\fi

\newsection Formulation of the Littlewood-Richardson rule.

\subsection Symmetric polynomials and partitions.

We fix some $n\in\N$, and let $\setof X_i:i\in\upto n\endset$ be a set of
$n$~indeterminates. The symmetric group~$\Sym_n$ acts on this set, and hence
on the ring~$\Z[\li(X0..n-1)]$, by permuting the indeterminates. A polynomial
$P\in\Z[\li(X0..n-1)]$ is called symmetric if it is fixed by every
$\pi\in\Sym_n$. We shall denote by~$\SFn$ the set of symmetric polynomials,
which is a subring of~$\Z[\li(X0..n-1)]$ (the ring of invariants for the
action of~$\Sym_n$); since the the action of~$\Sym_n$ preserves the natural
grading of~$\Z[\li(X0..n-1)]$ by total degree, $\SFn$~is a graded ring, and we
shall denote by $\SFn^d$ the set of homogeneous symmetric polynomials of
degree~$d$. Our interest will be in an explicit description of the
multiplicative structure of this ring, expressed on a particular $\Z$-basis,
that of the so-called Schur polynomials. But let us first consider some other
$\Z$-bases.

The simplest symmetric polynomials are those which are formed as the sum over
an orbit of a monomial in~$\Z[\li(X0..n-1)]$; we shall these symmetric
monomials. The monomials in~$\Z[\li(X0..n-1)]$ are of the form
$X_0^{\alpha_0}\cdots X_{n-1}^{\alpha_n-1}$, which we shall abbreviate to
$X^\alpha$, where $\alpha=(\li(\alpha0..n-1))\in\N^n$. We shall write
$|\alpha|=\deg X^\alpha=\sum_{i\in\upto n}\alpha_i$. In order to select a
specific representative within each orbit of monomials, we define a partial
ordering~`$\preceq$' on the set of monomials~$X^\alpha$; or equivalently on
the set~$\N^n$ of multi-exponents~$\alpha$. This ordering is generated by the
relations $X_i\succ X_{i+1}$ and the condition that $X^\alpha\prec X^\beta$
implies $MX^\alpha\prec MX^\beta$ for any monomial~$M$ (this implies that
monomials of distinct degrees are always incomparable); explicitly, one has
$X^\alpha\prec X^\beta$ if and only if
$\sum_{i<k}\alpha_i\leq\sum_{i<k}\beta_i$ for $k\in\upto n$, and
$|\alpha|=|\beta|$. Every $\Sym_n$-orbit of monomials contains a maximum
for~`$\preceq$' , which is a monomial $X^\\$ with $\fold\geq(\\0..n-1)\geq0$;
we shall call such a $\\\in\N^n$ a \newterm{partition} of~$d=|\\|$ into
$n$~parts, written $\\\in\Part_{d,n}$. For such~$\\$ we define the symmetric
monomial $m_\\(n)=\sum_{\alpha\in\Sym_n\cdot\\}X^\alpha$, where
$\Sym_n\cdot\\$ denotes the $\Sym_n$-orbit of~$\\$ in~$\N^n$. The symmetric
monomials $m_\\(n)$, for $\\\in\Part_{d,n}$, form a $\Z$-basis of $\SFn^d$. A
partition~$\\$ of~$d$ (without qualification, written $\\\in\Part_d$) is
defined as a sequence $(\\_i)_{i\in\N}$ of natural numbers (called parts) with
$\\_i\geq\\_{i+1}$ for all~$i$, such that $\\_m=0$ for some~$m\in\N$, and
$\sum_{i\in\upto m}\\_i=d$. We shall identify $\Part_{d,n}$ with the subset of
$\Part_d$ of partitions that have at most~$n$ non-zero parts; we also put
$\Part=\Union_{d\in\N}\Part_d$.

In the special case that $\\$ is the partition of $d\leq n$ into $n$~parts for
which all ($d$) non-zero parts are~$1$, the symmetric monomial $m_\\(n)$ is
called the $d$-th elementary symmetric polynomial, and written~$e_d(n)$. The
``fundamental theorem on symmetric functions'' states that the polynomials
$e_i(n)$, for $1\leq i\leq n$, generate~$\SFn$ as a ring, and are
algebraically independent; in other words, $\SFn$ is isomorphic as a graded
ring to $\Z[\li(Y1..n)]$ with $\deg Y_i=i$, the isomorphism sending $Y_i$
to~$e_i(n)$. We shall however not make use of this fact (we refer
to~\ref{Macdonald} for many more facts about symmetric polynomials and
functions). We may also express $e_d(n)$
as~$\sum_{0\leq\fold<(i0..d-1)<n}\fold{}(Xi_0..i_{d-1})$. If in this
expression we replace the strict inequalities between the indices by weak
ones, then we obtain another symmetric polynomial
$h_d(n)=\sum_{0\leq\fold\leq(i0..d-1)<n}\fold{}(Xi_0..i_{d-1})$, called the
$d$-th complete symmetric polynomial. The fact that this is indeed a symmetric
polynomial follows from the fact that it contains every monomial $X^\alpha$
with $|\alpha|=d$ exactly once, whence
$h_d(n)=\sum_{\\\in\Part_{d,n}}m_\\(n)$. Note that if we would replace only
some strict inequalities by weak ones, the result would not be a symmetric
polynomial. Like the $e_i(n)$, the $h_i(n)$ for $1\leq i\leq n$ form a
polynomial basis for~$\SFn$. We saw that $\setof
m_\\(n):\\\in\Part_{d,n}\endset$ is a $\Z$-basis of~$\SFn$; since monomials of
degree~$d$ in the generators $e_i(n)$ or $h_i(n)$ are naturally parametrised
by~$\Part_{d,n}$, each part~$i$ representing a factor $e_i(n)$~or~$h_i(n)$, we
have two more such bases parametrised by the same set.

\subsection Semistandard tableaux and Schur polynomials.

The Schur polynomials of degree~$d$ form yet another $\Z$-basis of~$\SFn$
parametrised by~$\Part_{d,n}$. Although their significance is not immediately
obvious from a purely ring-theoretic perspective, they are of fundamental
importance in many situations where the ring~$\SFn$ is encountered. For
instance, $\SFn$ occurs as the character ring of polynomial $\GL_n(\C)$
representations, and the Schur polynomials are the irreducible characters.
They can be defined as quotients of alternating polynomials, or be expressed
in terms of power sums ($m_d(n)$) using symmetric group characters; for our
purposes however, a purely combinatorial description will serve best, and it
is in that way that we shall define Schur polynomials. For those who wish some
motivation for this definition, we refer to places where it is respectively
deduced from an algebraic definition \ref{Macdonald, I~(5.12)}, from an
explicit construction of irreducible $\GL_n(\C)$ representations \ref{Fulton},
and even in an axiomatic approach \ref{Zelevinsky Hopf algebra}.

Our definition of Schur polynomials is rather similar to the description of
$e_d(n)$~and~$h_d(n)$ as the sum of a collection of monomials
$\fold{}(Xi_0..i_{d-1})$; indeed $e_i(n)$~and~$h_i(n)$ are instances of Schur
polynomials. The linear sequence of strict respectively weak inequalities
relating the indices $\li(i0..n-1)$ for $e_i(n)$~and~$h_i(n)$ are replaced for
general Schur polynomials by a more complicated mixture of strict and weak
inequalities. As we remarked above, this does not always give rise to a
symmetric polynomial; we shall see however that it will do so when the
inequalities follow a specific pattern associated to certain $d$-element
subsets of~$\N^2$ called diagrams. For the Schur polynomial~$s_\\(n)$ (with
$\\\in\Part_d$), this will be the \newterm{Young diagram}~$Y(\\)$ of~$\\$,
defined as $\setof(i,j)\in\N^2:j\in\upto{\\_i}\endset$; an index of summation,
ranging over~$\upto n$, is associated with each of its elements. These
elements will be drawn as squares in the plane, and for a fixed term of the
summation, the associated index will be written into it; the arrangement is
like that of matrix entries, the square~$(i,j)$ being in row~$i$ and
column~$j$. One obtains a sequence of left-justified rows of successive
lengths $\\_0,\\_1,\ldots$; for instance, the Young diagram of
$\\=(4,2,1,0,\ldots)$ is drawn as $\smallsquares\smallsquares\Young(,,,|,|)$.
Later, we shall use the notion of \newterm{diagonals}, defined as sets
$\setof(i,j)\in\N^2:j-i=k\endset$ for some constant~$k$, called the
\newterm{index} of the diagonal; the set of all diagonals is totally ordered
by their indices.

The terms of the summation that will define~$s_\\(n)$ are determined by fixing
all indices of summation, thus giving an assignment $T\:(i,j)\mapsto T_{i,j}$;
these indices are subject to the inequalities (whenever defined)
$T_{i,j}<T_{i+1,j}$ (strict increase down columns) and $T_{i,j}\leq T_{i,j+1}$
(weak increase along rows). An assignment~$T$ satisfying these conditions is
called a \newterm{semistandard Young tableau} of shape~$\\$ and entries
in~$\upto n$, and the set of all such is denoted by~$\Tab(\\,n)$. As an
example, for $T={\smallsquares\smallsquares\Young(0,0,1,1|1,2|2)}$ one
has~$T\in\Tab((4,2,1),3)$. Like in the case of $e_d(n)$~and~$h_d(n)$, each
term in the summation will be the term of factors $X_k$, one for each
occurrence of~$k$ as summation index. Therefore we define for $T\in\Tab(\\,n)$
its \newterm{weight}~$\alpha=\wt T\in\N^n$ to be such that $\prod_{(i,j)\in
Y(\\)}X_{T_{i,j}}=X^\alpha$, in other words $\alpha_k$ counts the occurrences
of the entry~$k$ in~$T$; in the example given one has $\wt T=(2,3,2)$. Then
for~$\\\in\Part_{d,n}$, the Schur polynomial~$s_\\(n)$ is defined by
$$
  s_\\(n)=\sum_{T\in\Tab(\\,n)}x^{\wt T}.  \label(\Schurdef)
$$
For instance, by enumerating~$\Tab((4,2,1),3)$ one finds that $s_{(4,2,1)}(3)$
is a symmetric polynomial with $15$~terms, which equals
$m_{(4,2,1)}(3)+m_{(3,3,1)}(3)+2m_{(3,2,2)}(3)$. The $\Sym_n$-invariance
of~$s_\\(n)$ is not at all evident from the definition, however. Although this
will follow from properties independently derived later, let us proof this key
fact right now.

\proclaim Proposition. \Schursym
The Schur polynomials are symmetric polynomials, i.e., $s_\\(n)\in\SFn^d$ for
$\\\in\Part_{d,n}$.

\proof
It will suffice to prove for any $k<n-1$ that $s_\\(n)$ is invariant under the
interchange of $X_k$~and~$X_{k+1}$; to this end we construct an involution of
the set $\Tab(\\,n)$, that realises an interchange of the components
$\alpha_k$~and~$\alpha_{k+1}$ of the weight~$\alpha$. We shall leave all
entries $T_{i,j}\notin\{k,k+1\}$ unchanged, as well as the entries
$k$~and~$k+1$ in any column of~$T$ that contains both of them. It is readily
checked that the set of squares containing the remaining entries (i.e.,
entries $k$~or~$k+1$ that as such are unique in their column) meets any given
row in a contiguous sequence of squares. Let the entries of that sequence be
$r$ times~$k$ followed by $s$ times~$k+1$ ($r$~or~$s$ might be~$0$); we
replace them by $s$ times~$k$ followed by $r$ times~$k+1$. The transformation
of~$T$ consists of performing this change independently for each row; clearly
the operation is an involution, and has the desired effect on~$\wt T$.
\QED

\heading "Remark"
This involution, which was introduced in~\ref{Bender Knuth}, is simple to
describe, but not the best one from a mathematical point of view, as it does
not give rise to an action of~$\Sym_n$ on~$\Tab(\\,n)$ (denoting interchange
of $i$~and~$i+1$ by~$s_i$, one may check that application of $s_0s_1s_0$ and
$s_1s_0s_1$ to the tableau~$T$ shown above gives different results). There is
another involution, the \foreign{rel\`evement plaxique} of
\ref{Lascoux Schutzenberger, \Sec4}, that does extend to an $\Sym_n$-action;
it is related to the coplactic operations discussed in~\Sec3 below.

Other properties of Schur polynomials are easy to establish. One has
$\wt{T}\preceq\\$ for all $T\in\Tab(\\,n)$, with equality for exactly one
such~$T$, namely the tableau with $T_{i,j}=i$ for all~$(i,j)\in Y(\\)$; this
tableau will be denoted by~$\Can\\$. This fact follows from the observation
that in any $T\in\Tab(\\,n)\setminus\{\Can\\\}$, one can find at least one
entry~$k+1$ that can be replaced by~$k$ (for some~$k<n-1$). Being a symmetric
polynomial, $s_\\(n)$ can be written as a sum of terms $m_\mu(n)$ for
$\mu\in\Part_{d,n}$ with integer coefficients. These coefficients are all
non-negative, they are zero unless $\mu\preceq\\$, and the coefficient
of~$m_\\(n)$ is equal to~$1$: the transition from the $s_\\(n)$ to the
$m_\\(n)$ is ``unitriangular'' with respect to~`$\preceq$'. Then the fact that
$\setof m_\\(n):\\\in\Part_{d,n}\endset$ is a $\Z$-basis of~$\SFn$, implies
that $\setof s_\\(n):\\\in\Part_{d,n}\endset$ is one as well.

We can now state the problem with which the Littlewood-Richardson rule is
concerned, as expressing multiplication in~$\SFn$ on the basis of the Schur
polynomials. In somewhat more detail: given $\\\in\Part_{d,n}$ and
$\mu\in\Part_{d',n}$, we wish to determine the integer coefficients
$c_{\\,\mu}^\nu$ for all~$\nu\in\Part_{d+d',n}$, such that
$$
  s_\\(n)s_\mu(n)=\sum_{\nu\in\Part_{d+d',n}}c_{\\,\mu}^\nu s_\nu(n).
  \label(\LRcdef)
$$
We have suppressed~$n$ in the notation $c_{\\,\mu}^\nu$, since it will turn
out that this coefficient is independent of~$n$ (although $n$ must be
sufficiently large for $c_{\\,\mu}^\nu$ to appear in the formula in the first
place). More generally, for identities valid for any number~$n$ of
indeterminates, we shall sometimes write $h_d$ for~$h_d(n)$ and $s_\\$
for~$s_\\(n)$, etc.\ (there is an algebraic structure called the ring of
symmetric functions that justifies this notation, see~\ref{Macdonald}, but we
shall not discuss it here). The~$c_{\\,\mu}^\nu$ are called
Littlewood-Richardson coefficients. Representation theoretic considerations
show that $c_{\\,\mu}^\nu\in\N$; the Littlewood-Richardson rule will in fact
describe the $c_{\\,\mu}^\nu$ as the cardinalities of certain
combinatorially defined sets.

\subsection Skew shapes, and skew Schur polynomials.

In order to formulate the Littlewood-Richardson rule, we need to extend the
class of tableaux beyond that of the semistandard Young tableaux.
Giving~$\N^2$ its natural partial ordering (simultaneous comparison of both
coordinates), Young diagrams can be characterised as its finite order ideals
(if $r\in Y(\\)$, then also $q\in Y(\\)$ for all $q\leq r$). For the class of
skew diagrams, this is relaxed to convexity with respect to~`$\leq$': a
\newterm{skew diagram}~$D$ is a finite subset of~$\N^2$ for which $p,r\in D$
and $p\leq q\leq r$ imply~$q\in D$. The summation analogous to (\Schurdef)
using a skew diagram will still give rise to a symmetric polynomial (one may
for instance check that the proof of proposition~\Schursym\ remains valid). A
skew diagram can always be written as the difference of two Young diagrams:
$D=Y(\\)\setminus Y(\mu)$ with $Y(\mu)\subset Y(\\)$; this representation is
not unique in general (although in some cases it is, for instance when the
sets of rows and columns meeting~$D$ are both initial intervals of~$\N$). In
many cases, for instance when considering the shapes of tableaux, it will be
important to fix the partitions $\\,\mu$ used to represent a skew diagram,
which leads us to define the related but distinct notion of a skew shape.

A \newterm{skew shape} is a symbol $\\/\mu$ where $\\,\mu\in\Part$ and
$Y(\mu)\subset Y(\\)$; the set of all skew shapes will be denoted by~$\Sk$.
For $\\/\mu\in\Sk$ we define $\Ylm=Y(\\)\setminus Y(\mu)$ and
$|\\/\mu|=|\Ylm|=|\\|-|\mu|$. A shape of the form~$\\/0$ is called a partition
shape; occasionally we consider partitions as skew shapes, in which case $\\$
is identified with~$\\/0$. It is a trivial but useful fact that if
$|\\/\mu|\leq3$, no diagonal can meet $\Ylm$ in more than one square. We shall
say that a skew shape~$\chi$ represents the product $\chi_0*\chi_1$ of two
other skew shapes if $Y(\chi)$ can be written as a disjoint union of skew
diagrams $\bar\chi_0$ and~$\bar\chi_1$, where $\bar\chi_k$ is obtained by some
translation from~$Y(\chi_k)$ ($k=0,1$), while for all $(i,j)\in\bar\chi_0$ and
$(i',j')\in \bar\chi_1$ one has $i>i'$ and $j<j'$. For instance,
$$ {\smallsquares \Skew(7:|5:,,|5:|4:|0:,|0:,|0:)}
\text{represents}
 {\smallsquares \Young(,|,|)}\quad*\quad{\smallsquares \Skew(3:|1:,,|1:|0:)}
$$
One could define an equivalence relation on~$\Sk$ such that this relation
defines a monoid structure on the quotient set, but this would be cumbersome,
and is not really needed for our purposes.

Let $\chi=\\/\mu\in\Sk$; a \newterm{skew semistandard tableau}~$T$ of
shape~$\chi$ and with entries in~$\upto n$ is given by specifying~$\chi$
itself, together with a map $Y(\chi)\to\N$ written $(i,j)\mapsto T_{i,j}$,
satisfying $T_{i,j}\in\upto n$, \ $T_{i,j}<T_{i+1,j}$ and $T_{i,j}\leq
T_{i,j+1}$ whenever these values are defined. The set of all such~$T$ is
denoted by~$\Tab(\chi,n)$, and $\Tab(\chi)=\Union_{n\in\N}\Tab(\chi,n)$; we
identify $\Tab(\\)$ with~$\Tab(\\/0)$. The weight $\alpha=\wt T\in\N^n$ of
$T\in\Tab(\chi,n)$ is defined by $\prod_{(i,j)\in
Y(\chi)}X_{T_{i,j}}=X^\alpha$, and the skew Schur polynomial $s_\chi(n)$ by
$$
  s_\chi(n)=\sum_{T\in\Tab(\chi,n)}x^{\wt T}.  \label(\skewSchurdef)
$$
Note that unlike (\Schurdef), there is no restriction to~$\Part_{d,n}$ here,
for the partitions $\\,\mu$ forming~$\chi$; the skew Schur polynomial will be
non-zero as long as there are no columns in~$Y(\chi)$ of length exceeding~$n$.
If a skew shape~$\chi$ represents the product $\chi_0*\chi_1$ of two other
skew shapes, then there is an obvious weight preserving bijection
$\Tab(\chi)\to\Tab(\chi_0)\times\Tab(\chi_1)$; therefore
$s_\chi=s_{\chi_0}s_{\chi_1}$. In this case we shall shall denote~$s_\chi$
by~$s_{\chi_0*\chi_1}$, or if $\chi_0=(\\/0)$ and $\chi_1=(\mu/0)$,
by~$s_{\\*\mu}$; our problem can be restated as finding the decomposition of
the skew Schur polynomial $s_{\\*\mu}(n)$ as a sum of ordinary Schur
polynomials. In fact the Littlewood-Richardson rule will describe the
decomposition of any skew Schur polynomial $s_\chi(n)$. We define coefficients
~$c_\chi^\nu$ for any skew shape~$\chi$ and $\nu\in\Part_{|\chi|}$, by the
decomposition formula
$$
  s_\chi=\sum_{\nu\in\Part_{|\chi|}}c_\chi^\nu s_\nu.
  \label(\skewLRcdef)
$$
If $\chi$ represents $(\\/0)*(\mu/0)$, then we have
$c_\chi^\nu=c_{\\,\mu}^\nu$.

\subsection Littlewood-Richardson tableaux.

Let $\\/\mu\in\Sk$, and $\nu\in\Part$ with $|\\/\mu|=|\nu|$; we shall now
introduce the objects that are counted by~$c_{\\/\mu}^\nu$, which will be
called Littlewood-Richardson tableaux of shape~$\\/\mu$ and weight~$\nu$. We
need a preliminary definition. For $T\in\Tab(\\/\mu)$ and $k,l\in\N$, define
$T_k^l$ to be the number of entries~$l$ in row~$k$ of~$T$, i.e., the
cardinality of the set $\setof j\in\upto{\\_k}\setminus\upto{\mu_k}:
T_{k,j}=l\endset$.

\proclaim Definition. \companiondef
Two tableaux $T\in\Tab(\\/\mu)$ and~$\bT\in\Tab(\nu/\kappa)$ are called
companion tableaux if $T_k^l=\bT_l^k$ for every $k,l\in\N$. In this case
$T$ is called $\nu/\kappa$-dominant (and $\bT$ is $\\/\mu$-dominant).

Here is an example of a pair of companion tableaux, with a table of the
pertinent values $T_k^l=\bT_l^k$:
$$
  T={\maybesmall\Skew(3:0,1|1:0,1,1,3|0:0,2,2,3|0:1,4,4,5|0:3,5)},
\quad
  \bT={\maybesmall\Skew(6:0,1,2|4:0,1,1,3|4:2,2|2:1,2,4|1:3,3|0:3,4)},
\quad
\(T_k^l\)_{\scriptstyle 0\leq k<5\atop\scriptstyle0\leq l<6}=
 \pmatrix{1&1&0&0&0&0\cr
          1&2&0&1&0&0\cr
	  1&0&2&1&0&0\cr
	  0&1&0&0&2&1\cr
	  0&0&0&1&0&1\cr}
 \label(\companionex)
$$
Note that if $T$ is $\nu/\kappa$-dominant, then $\wt T=\nu-\kappa$, since
$(\wt{T})_j=\sum_iT_i^j=\sum_i\bT_j^i=\nu_j-\kappa_j$. We shall simply say
that $T$ is $\kappa$-dominant if it is $\nu/\kappa$-dominant for
$\nu=\kappa+\wt T$; this notion is used in~\ref{Littelmann paths}. $T$~may
have companion tableaux of different shapes, but at most one of a given
shape~$\nu/\kappa$: the multiset of entries of any row is determined by~$T$,
and they must be weakly increasing. Therefore different companion tableaux
of~$T$ differ only by horizontal slides of their rows; in the above example
one could for instance shift the first~$4$ rows of~$\bT$ one place to the
left. To test for $\kappa$-dominance of a tableau~$T$, it suffices to
construct the unique candidate for~$\bT$ and check that its columns are
strictly increasing.

\proclaim Definition. \LRtabdef
A tableau $L\in\Tab(\\/\mu)$ is a Littlewood-Richardson tableau if it is
$0$-dominant. The set of all $\nu/0$-dominant tableaux in~$\Tab(\\/\mu)$ is
denoted by~$\LRt\\/\mu,\nu$.

Here is an example Littlewood-Richardson tableau~$L$, and its companion Young
tableau~$\bar L$.
$$
  L=\Skew(4:0,0|2:0,1,1|1:0,1,2,2|0:0,1,3|0:2,4),
\qquad
  \bar L=\Young(0,0,1,2,3|1,1,2,3|2,2,4|3|4).
  \label(\LRtablex)
$$
The given definition of Littlewood-Richardson tableaux is not the traditional
one, but that definition can be easily derived from it. A pair of companion
tableaux $T\in\Tab(\\/\mu)$ and~$\bT\in\Tab(\nu/\kappa)$ determines a
bijection $p\:\Ylm\to Y(\nu/\kappa)$, if one fixes for each~$k,l$ a bijection
between the sets of squares whose entries are counted by~$T_k^l$ and
$\bT_l^k$. Both these sets are contiguous sequences of squares of some row; we
choose the bijection that \emph{reverses} the left-to-right ordering of the
squares. For the companion tableaux $T,\bT$ in~(\companionex) we get the
following bijection, indicated by matching labels:
$$
  p\: \Skew(3:b,a|1:f,e,d,c|0:l,k,h,g|0:q,p,n,m|0:s,r)
\longrightarrow
  \Skew(6:b,f,l|4:a,d,e,q|4:h,k|2:c,g,s|1:n,p|0:m,r)
 \label(\pictureex)
$$
We have the following property: for any square~$(i,j)\in\Ylm$, all squares
whose images are on the same row as~$p(i,j)$ and to the left of it, themselves
lie in a column to the right of~$(i,j)$, and in the same row or above it; in
formula, if $p(i,j)=(r,c)$ and $p(i',j')=(r,c')$ with $c'<c$, then $j'>j$ and
$i'\leq i$. E.g., for the square~$(1,2)$ labelled~$e$ on the left in the
example, the labels $d$~and~$a$ appear strictly to its right and weakly above
it. Given $T$ and~$\kappa$, one can therefore construct~$p$---and
implicitly~$\bT$---by traversing~$\Ylm$ in such an order that the mentioned
other squares~$(i',j')$ are always encountered before~$(i,j)$ is (in the
example the alphabetic order of the labels has this property): the
image~$p(i,j)$ can then be taken to be the first currently unused square of
row~$T_{i,j}$ of~$Y(\nu/\kappa)$. If moreover the traversal is by rows (like
the alphabetic order in the example), then strict increase in the column
of~$\bT$ at the square~$p(i,j)$ can be checked as soon as $p(i,j)$ is located:
if the square directly above~$p(i,j)$ lies in~$Y(\nu/\kappa)$, then
it must have been included in the image of~$p$ before~$p(i,j)$ is. What this
amounts to, is that at each point during the construction, the union of
$Y(\kappa)$ and the image of~$p$ so far constructed must be a Young diagram.

\proclaim Proposition. \LRlatticechar
A tableau~$T\in\Tab(\\/\mu,n)$ is $\kappa$-dominant if and only if the
following test succeeds. A variable~$\alpha\in\N^n$ is initialised
to~$\kappa$; then the squares~$(i,j)\in\Ylm$ are traversed by weakly
increasing~$i$, and for fixed~$i$ by strictly decreasing~$j$: at
square~$(i,j)$ the component $\alpha_{T_{i,j}}$ is increased by~$1$. The test
succeeds if and only if one has $\alpha\in\Part$ throughout the entire
procedure.
\QED

For~$\kappa=0$ this is still not quite the traditional description of
Littlewood-Richardson tableaux, which states that the word over the alphabet
$\upto n$ obtained by listing the entries~$T_{i,j}$ in the order described in
the proposition (e.g., for the Littlewood-Richardson tableau~$L$ shown
in~(\LRtablex), the word $\let~=\,0~0~1~1~0~2~2~1~0~3~1~0~4~2$) should be a
``lattice permutation'' (see \Sec3 for a definition; one also finds the terms
``lattice word'' and ``Yamanouchi word''). However, the characterisation above
results directly from expanding the definition of that term. Incidentally, the
(very old) term lattice permutation appears to be related to the fact that the
sequence of values assumed by~$\alpha$ in the proposition describes a path
(from~$\kappa$ to~$\kappa+\wt T$) in the lattice~$\N^n$, that is confined to
remain inside the cone~$\Union_{d\in\N}\Part_{d,n}\subset\N^n$. It may also be
noted that the \emph{original} formulation in \ref{Littlewood Richardson} is
very close to our definition~\LRtabdef, see~\Sec4. We can now state the Rule.

\proclaim Theorem [Littlewood-Richardson rule]. \LRrule
For all $\chi\in\Sk$ and $\nu\in\Part$ one has
$c_\chi^\nu=\Card{\LRt\chi,\nu}$.

We shall present a proof of this theorem in the course of \Sec\Sec$2,3$. In
view of the definitions (\skewLRcdef) of~$c_\chi$ and
(\Schurdef)~and~(\skewSchurdef) of (skew) Schur polynomials, it suffices to
construct a bijection
$$
  \Rob\:\Tab(\chi,n)\to
  \Union_{\nu\in\Part_{|\chi|,n}}\LRt\chi,\nu\times\Tab(\nu,n)
  \label(\Robinson)
$$
for any~$n$, such that if $\Rob(T)=(L,P)$ one has $\wt T=\wt P$. It will then
follow from the fact that the same set $\LRt\chi,\nu$ occurs, independently
of~$n$ (provided only that $\nu$ occurs in the summation), that the
coefficient~$c_\chi^\nu$ is independent of~$n$, as we claimed. We shall
refer to~$\Rob$ as Robinson's correspondence, since it was first described
in~\ref{Robinson} (albeit in different terms, not using tableaux). It is
natural to describe the correspondence by separately defining its components
$\Rob_0$ and~$\Rob_1$, where $\Rob(T)=(\Rob_0(T),\Rob_1(T))$.

The decomposition of ordinary Schur functions is trivial, so according to the
Littlewood-Richardson rule, $\LRt\\/0,\nu$ should be empty unless $\\=\nu$, in
which case it should be a singleton. Indeed, for companion tableaux
$T\in\Tab(\\/0)$ and $\bT\in\Tab(\nu/0)$ one has $\wt T=\nu$ and $\wt\bT=\\$,
which together with $\wt T\preceq\\$ and~$\wt\bT\preceq\nu$ imply $\\=\nu$ and
$T=\bT=\Can\\$. If $L\in\LRt\chi,\nu$ is a Littlewood-Richardson tableau whose
shape~$\chi$ represents~$\\*\mu$, then the test of proposition~\LRlatticechar\
for $0$-dominance shows that the restriction~$L_\mu$ of~$L$ to the
squares belonging to the factor~$\mu$ of~$\\*\mu$ is itself a
Littlewood-Richardson tableau, and hence equal (up to translation)
to~$\Can\mu$; then the restriction~$L_\\$ of~$L$ to the remaining squares is
$\nu/\mu$-dominant. By the symmetry of companion tableaux, this proves:

\proclaim Proposition. \LRswitchprop
$\Card{\LRt\chi,\nu}=\Card{\LRt\nu/\mu,\\}$ for~$\chi\in\Sk$
representing $\\*\mu$.
\QED

As an example of the correspondence underlying this propostiion, we derive
from~(\LRtablex) the following pair of corresponding Littlewood-Richardson
tableaux, for $\\=(5,4,3,1,1)$, $\mu=(4,2,1)$, $\nu=(6,5,5,3,2)$:
$$
  \Skew(5:0,0,0,0|5:1,1|5:2|0:0,0,1,2,3|0:1,1,2,3|0:2,2,4|0:3|0:4)
\longleftrightarrow\qquad
  \Skew(4:0,0|2:0,1,1|1:0,1,2,2|0:0,1,3|0:2,4)
$$

In order to use the Littlewood-Richardson rule to decompose a skew Schur
polynomial~$s_\chi$ on the basis of Schur polynomials, an effective
enumeration is  required of the union of sets~$\LRt\chi,\nu$, as $\nu$ varies
over~$\Part_{|\chi|}$. This can be done using an efficient search procedure,
in which the choices for the entries of the squares of~$Y(\chi)$ are fixed in
the same order as the traversal of proposition~\LRlatticechar\ (alternatively,
any ``valid reading order'' as defined below works equally well). For any
square~$(i,j)$, the possible values~$T_{i,j}$ that can be chosen must make the
test of that proposition succeed (i.e., they must index a part of the
partition~$\alpha$ constructed so far that can be increased), and they must
satisfy the monotonicity conditions for the rows and columns of~$T$. The
reason that we called the search efficient, is that there is always at least
one value that satisfies all these conditions (cf.\ \ref{pictures,
proposition~2.5.3}), so that the search tree will not have unproductive
branches. Here we have assumed that no upper limit~$n$ is imposed on the
number of non-zero parts of~$\nu$ and hence on the entries in the
Littlewood-Richardson tableau~$T$; such a restriction can however be
incorporated into the search, by placing an additional condition on the
values~$T_{i,j}$ tried: they should be sufficiently small to allow column~$j$
of~$T$ to be completed in a strictly increasing manner. With this extra
requirement there will still always remain at least one possible value,
unless~$Y(\chi)$ has columns of length exceeding~$n$ (in which case of course
no suitable tableaux~$T$ exist at~all).

By contrast, no efficient search procedure is known for enumerating just a
single set~$\LRt\chi,\nu$, i.e., one for which the size of the search tree
that has to be traversed is proportional to~$\Card{\LRt\chi,\nu}$. In fact no
sufficient condition for $\LRt\chi,\nu$ to be non-empty, that does not amount
to actually finding an element, is even known; therefore it is not possible in
general to tell beforehand whether some value tried for an entry will lead to
any solutions. We observe that, while all bijections described in this paper
are readily computed, it is in some cases much harder to construct the sets
themselves linked by these bijections. We add one more remark, prompted by the
fact that many texts give the Littlewood-Richardson rule only in the form
$c_{\\,\mu}^\nu=\Card{\LRt\nu/\mu,\\}$, while defining Littlewood-Richardson
tableaux in terms of a reading of the tableau. This would falsely suggest that
the rule is not practical for calculating a product~$s_\\s_\mu$, by requiring
either the construction of complete trial tableaux of varying shapes before
testing the Littlewood-Richardson condition (which would fail in most cases),
or separate searches for any plausible shape~$\nu/\mu$ and fixed weight~$\\$
(which is also unattractive for reasons just mentioned). In reality, viewing
the search strategy outlined above for Littlewood-Richardson tableaux of a
shape representing~$\\*\mu$ from the perspective of the companion tableaux,
one finds an efficient enumeration procedure for $\Union_\nu\LRt\nu/\mu,\\$.
Actually, such an enumeration procedure is just what the rule, in its original
form given by Littlewood and Richardson, describes; a literal quotation of
this description can be found in~\Sec4.

\subsection Pictures and reading orders.
\labelsec \readingsec

We have constructed a
bijection $p\:\Ylm\to Y(\nu/\kappa)$ corresponding to a pair of companion
tableaux $T\in\Tab(\\/\mu)$ and~$\bT\in\Tab(\nu/\kappa)$, in order to derive proposition~\LRlatticechar. Since $T$~and~$\bT$
can easily be reconstructed from~$p$, one may study the bijections that arise
in this way, in place of such pairs of companion tableaux, or of
$\nu/\kappa$-dominant tableaux of shape~$\\/\mu$. The conditions that
$T$~and~$\bT$ be semistandard tableaux translate into a geometric
characterisation of these bijections; this leads to the concept
of~\newterm{pictures} introduced in \ref{Zelevinsky pictures}. Due to the
symmetry between companion tableaux, the inverse of any picture is again a
picture. One important aspect of pictures is that there are many equivalent
ways to define them (like the different characterisations of companion
tableaux above). We shall not discuss pictures in depth here, for which we
refer to \ref{Fomin Greene} and~\ref{pictures}, but it is useful to make a few
observations that result from study of pictures.

We have observed above that if one traverses a row from left to right, then
the row index of the (inverse) image by~$p$ increases weakly, while the column
index decreases strictly. A similar condition holds for traversal of a column
from top to bottom: the row index of the (inverse) image increases strictly,
while the column index decreases weakly (this can be proved by induction on
the column considered). Therefore, when verifying $\kappa$-dominance with a
single traversal of~$T$ (as in proposition~\LRlatticechar), it is not
necessary to have encountered all squares in rows above it before handling a
square~$(i,j)$: only those in column~$j$ or to its right can influence the
test made at~$(i,j)$. We might for instance also traverse columns from top to
bottom, processing the columns from right to left. We shall consider any
ordering of the squares that encounters $(i',j')$ before~$(i,j)$ whenever
$i'\leq i$ and $j'\geq j$ to be a \newterm{valid reading order}. The two cases
where one proceeds systematically by rows or by columns merit special names:
we shall refer to the former as the \newterm{Semitic reading order} (after the
Arabic and Hebrew way of writing), and to the latter as the \newterm{Kanji
reading order} (after the Japanese word for Chinese characters, which are
thusly read). In~(\pictureex), the Semitic reading of the left diagram gives
$ab cdef ghkl mnpq rs$, the Kanji reading order gives $ac bdgm ehn fkpr lqs$,
while $acbdegfhkmlnprqs$ corresponds to yet another a valid reading order.

In constructing the picture~$p$ corresponding to a $\nu/\kappa$-dominant
tableau~$T$, a distinction is forced between the values at all squares
of~$\Ylm$, i.e., $p$ is injective even if~$T$ is not. Such a distinction can
be used in order to apply to~$T$ operations that are initially defined only
for tableaux with distinct entries (\Sec2 provides an example): it suffices to
use an injective map $r:Y(\nu/\kappa)\to\N$ such that $r\after p$ is a
tableau. Taking for~$r$ the map corresponding to the Semitic reading of
$Y(\nu/\kappa)$, one obtains a tableau $S=r\after p$ such that
$T_{i,j}<T_{i',j'}$ implies $S_{i,j}<S_{i',j'}$, and when $T_{i,j}=T_{i',j'}$
one has $S_{i,j}<S_{i',j'}$ if and only if $j<j'$; this is essentially the
operation of standardisation defined in~\Sec2. We note however that $r\after
p$ will be a tableau whenever $r$ corresponds to \emph{any} valid
reading~order; such tableaux may be called specialisations of~$p$. It is shown
in~\ref{pictures} that all relevant operations that are defined for~$T$ in
terms of its standardisation~$S$ could equally well be defined in terms of any
specialisation of~$p$. When studying semistandard tableaux however, it is
simplest to use the standardisation (as we shall do), mainly because it does
not depend on any choice of a shape $\nu/\kappa$ for which $T$ is
$\nu/\kappa$-dominant.

One can develop the theory of the Littlewood-Richardson rule entirely in terms
of pictures; doing so clarifies the structure behind many operations and makes
certain symmetries explicit. Nevertheless we believe the exposition in terms
of tableaux that we shall give is easier to understand, which is in part due
to the fact that pictures are (despite their name) more difficult to visualise
then tableaux are.

\endofchapter

 \par

\ifx\macrosloaded\undefined
  \input marxmax
  

\newif\ifnarrow
\ifdim \hsize<14.5cm \narrowtrue \else\narrowfalse \fi

\def\maybesmall{\ifnarrow\smallsquares\fi}

\def\ol#1
{\vbox{\vskip0.3pt \hrule \vskip 0.8pt \hbox{$#1$}}}

\def\ul#1
{\vtop{\hbox{$#1$}\vskip 0.8pt \hrule \vskip0.3pt }}

\font\frak=eufm10
\font\csc=cmcsc10
\font\URLfont=cmtex10 
\hyphenchar\URLfont=-1 

\mathorddef\Part=P2 
\mathorddef\A=A0 \mathorddef\B=B0 
\def\GL{\boldmath{GL}}
\def\Ylat{(\Part,{\subset})}
\def\upto#1{[#1]}
\opdef\Pic
\altopdef\Tab{SST}
\opdef\ST
\opdef\LR
\opdef\JdT
\opdef\wt
\mathorddef\T=T2 
\mathorddef\Sk=S2 
\def\SFn{\Lambda_n}
\def\Can#1{\boldmath1_{#1}}
\def\Ylm{Y(\\/\mu)}
\def\bT{\bar T}
\def\LRt#1,#2{\LR(#1,#2)}
\def\JdTt#1,#2{\Tab(#1)^{\slid{#2}}}
\mathorddef\Rob=R2 
\def\StT#1{\ST(#1)}
\def\snu{^{\slid\nu}}
\def\Uqgln{U_q(\hbox{\frak gl}_n)}
\def\leqr{\leq_r}
\def\leqKan{\leq_{\romath K}}

\def\wKan{w_{\romath K}}
\def\wSem{w_{\romath S}}

\mathdef\slid=Rel \triangleright
\mathdef\dils=Rel \triangleleft

\long\def\quotation#1
{{\smallskip\narrower\it\noindent\llap{``\/}#1\unskip''\smallskip}\noindent}

{\catcode`\.=\active\catcode`\/=\active\catcode`\_=\active
 \gdef\URL
 {\begingroup\def~{\char`\~}%
  \catcode`\.=\active\catcode`\/=\active\catcode`\_=\active
  \let_=\_%
  \def.{\char`\.\penalty\hyphenpenalty}\def/{\char`\/\penalty\hyphenpenalty}%
  \URLfont\URLnext
}}
\def\URLnext#1{http://#1\endgroup}

\hyphenation
{semi-standard Yamano-uchi mono-mial mono-mials inter-related combi-na-torics
 combi-na-torial algo-rithm valid-ity}

  \writelab{}
  \def\jobname{lrr}
\readrefs
  \let\endofchapter=\bye
  \input lrr1.lab
\fi

\newsection Tableau switching and jeu de taquin.

In this section we shall consider a combinatorial procedure that will turn out
to be intimately related to the Littlewood-Richardson rule. This procedure is
essentially Sch\"utzenbergers \foreign{jeu de taquin}, but we prefer to
introduce it in a slightly different form called ``tableau switching'' (a term
that was introduced in~\ref{tableau switching} for an operation that, although
defined in a slightly different way, constructs the same correspondence
between pairs of skew tableaux as we shall do below).

\subsection Chains in the Young lattice and standardisation.

Inclusion of Young diagrams defines a partial ordering on the set~$\Part$ of
partitions, which shall be denoted by~`$\subset$'; the poset~$\Ylat$ is called
the Young lattice. We define a \newterm{skew standard tableau} of
shape~$\\/\mu$ to be a saturated increasing chain in the Young lattice
from~$\mu$ to~$\\$, i.e., a sequence of partitions starting with~$\mu$ and
ending with~$\\$ such that the Young diagram of each partition is obtained
from that of its predecessor by the addition of exactly one square. We shall
denote the set of all skew standard tableaux of shape~$\\/\mu$
by~$\StT{\\/\mu}$, and say that each of its elements has size~$|\\/\mu|$. A
skew standard tableau of partition shape is called a standard Young tableau,
and we write $\ST(\\)$ for~$\ST(\\/0)$. If its shape~$\\/\mu$ is given, then
specifying some $S\in\StT{\\/\mu}$ amounts to putting a total ordering
on~$\Ylm$, describing the order in which the squares are added. This can be
done by labelling the squares in the desired order with increasing numbers (or
elements of some other totally ordered set); these labels will increase along
each row and column, whence the name tableau. For instance, the sequence with
Young diagrams 
$\smallsquares\smallsquares
\Young(,,|),~
\Young(,,|,),~
\Young(,,|,,),~
\Young(,,|,,|),~
\Young(,,,|,,|),~
\Young(,,,|,,|,)
$ is represented by $\smallsquares\smallsquares\Skew(3:3|1:0,1|0:2,4)$. In the
literature it is usually this representation that is called a standard
tableau, but chains of partitions will be more convenient for us to work with.

In fact any skew semistandard tableau~$T$ determines a chain of partitions,
with the convention that the squares~$(i,j)$ are ordered by increasing value
of their entries~$T_{i,j}$, or in case these entries are equal, by increasing
column number~$j$; this chain (a skew standard tableau) will be called the
standardisation of~$T$. For instance, the standardisation of
${\smallsquares\smallsquares\Skew(3:2|1:0,0|0:2,3)}$ is the skew standard
tableau depicted above. Any skew semistandard tableau is determined by its
standardisation together with its weight, but a given combination of a skew
standard tableau~$S$ and a weight~$\alpha$ does not necessarily correspond to
any skew semistandard tableau; if it does, we shall say that $S$ is compatible
with~$\alpha$.

\subsection Tableau switching and jeu de taquin.

Among the skew shapes~$\\/\mu$ of with $|\\/\mu|=2$, two different kinds can
be distinguished: shapes for which the two squares of $\Ylm$ are incomparable
in the natural ordering of~$\N^2$, and which shapes therefore admit two
different skew standard tableaux, and shapes for which those squares lie in
the same row or column (and are adjacent), which shapes admit only one skew
standard tableau. The latter kind of shapes will be called dominos. We shall
give a construction, based on a certain class of doubly indexed families of
partitions, that can be found in~\ref{eljc, 2.1.2}. Let $I=\{k,\ldots,l\}$ and
$J=\{m,\ldots,n\}$ be intervals in~$\Z$, and let $(\\^{[i,j]})_{i\in I,j\in
J}$ be a family of partitions; we shall call this a tableau switching family
on~$I\times J$ if each ``row'' $\\^{[i,m]},\ldots,\\^{[i,n]}$ and each
``column'' $\\^{[k,j]},\ldots,\\^{[l,j]}$ is a skew standard tableau, and if
$\\^{[i,j+1]}\neq\\^{[i+1,j]}$ whenever $\\^{[i+1,j+1]}/\\^{[i,j]}$ is not a
domino. Here is a small example:
$$
\(\\^{[i,j]}\)_{0\leq i,j<4}=\left(\enspace
{\smallsquares\smallsquares
  \dimen0=\ht\strutbox \advance\dimen0 by \dp\strutbox
  \vcenter{\halign
   {     $\vcenter to 18pt{\Young#\vss}$&&
    \quad$\vcenter to 18pt{\Young#\vss}$\cr
    \omit $\vbox to 9pt{\hbox{$\scriptscriptstyle\circ$}\vss}$
          & ()    & (,)    & (,|) \cr
     ()   & (|)   & (,|)   & (,|,) \cr
     (,)  & (,|)  & (,,|)  & (,,|,) \cr
     (,|) & (,||) & (,,||) & (,,|,|) \cr
   }}
}
\enspace\right)
 \label(\tabswitchex)
$$
Whenever either the sequence $\\^{[i,j]}$, $\\^{[i,j+1]}$, $\\^{[i+1,j+1]}$ or
the sequence $\\^{[i,j]}$, $\\^{[i+1,j]}$, $\\^{[i+1,j+1]}$ is specified to be
some skew standard tableau of size~$2$, there is a unique value for the
remaining partition for which one obtains a tableau switching family
on~$\{i,i+1\}\times\{j,j+1\}$ (because if $\\^{[i+1,j+1]}/\\^{[i,j]}$ is a
domino, then necessarily $\\^{[i,j+1]}=\\^{[i+1,j]}$). It follows that if
$\\^{[i,j]}$ is prescribed on some ``lattice path'' from $[k,m]$~to~$[l,n]$
(i.e., at each step one of the indices increases by~$1$) by any skew standard
tableau, then there is a unique way to extend this to tableau switching family
on~$I\times J$. This allows us in particular to make the following definition.

\proclaim Definition. \tabswitchdef
Let $S$ and $T$ be skew standard tableaux of respective shapes $\mu/\nu$
and~$\\/\mu$. The pair $(T',S')$ of skew standard tableaux obtained
from~$(S,T)$ by tableau switching, written $(T',S')=X(S,T)$, is defined by the
existence of a tableau switching family
$(\\^{[i,j]})_{k\leq i\leq l;~m\leq j\leq n}$ such that
$$
\eqalign
{S&=(\\^{[k,m]},\dots,\\^{[l,m]}),\cr
 T'&=(\\^{[k,m]},\dots,\\^{[k,n]}),\cr}
\quad
\eqalign
{T &=(\\^{[l,m]},\dots,\\^{[l,n]});\cr
 S'&=(\\^{[k,n]},\dots,\\^{[l,n]}).\cr}
$$

As an example, the tableau switching family~(\tabswitchex) establishes the
fact that
$X({\smallsquares\smallsquares\Young(0,1|2)}
  ,{\smallsquares\smallsquares\Skew(2:1|1:2|0:0)}
  )=
  ({\smallsquares\smallsquares\Young(0,1|2)}
  ,{\smallsquares\smallsquares\Skew(2:1|1:0|0:2)
  })$.
A visual way to interpret this definition is the following. Represent $S$ and
$T$ by labelling the squares of their diagrams, using two disjoint totally
ordered sets $A=\{\fold<(ak..l-1)\},B=\{\fold<(bm..n-1)\}$ of labels, say red
and blue numbers. Associating the~$a_i$ to vertical path segments and
the~$b_j$ to horizontal ones, each lattice path from $[k,m]$~to~$[l,n]$
corresponds to a shuffle of the sets $A$~and~$B$, i.e., a total ordering on
$A\union B$ that is compatible with the orderings of $A$~and~$B$ individually.
To each lattice path is associated a labelling of the squares of $Y(\\/\nu)$
with $A\union B$, that represents, for the corresponding total ordering of
$A\union B$, the skew standard tableau consisting of the $\\^{[i,j]}$ read off
along the path. According to the definition, one must move from the lattice
path via $[l,m]$ to the one via~$[k,n]$, while updating the labelling: the
initial labelling corresponds to $S$~and~$T$, and the final one to
$T'$~and~$S'$. Changing the lattice path in one place means transposing
one~$a_i$ with one~$b_j$ in the shuffle; according to
definition~\tabswitchdef, if the squares labelled by $a_i$~and~$b_j$ are
adjacent then their labels are interchanged, and otherwise the labelling
remains unchanged. In the former case the exchange of labels is a special case
of a switch in the sense of \ref{tableau switching}; this shows that $X(S,T)$
coincides with the result of the tableau switching procedure defined there (in
fact that procedure is more liberal, and allows using intermediate labellings
that do not correspond to a lattice path). We display a (not very systematic)
sequence of shuffles with the corresponding tableaux for the
example~(\tabswitchex), using italics for red numbers and bold face for blue
ones:
$$\def\r#1{{\it#1}}\def\b#1{{\bf#1}}
\everycr{\noalign{\hss}}
\maybesmall
\line
{\valign{\hbox{$\strut#$}\tabskip=4pt&\Young#\tabskip=0pt\cr
\r0\r1\r2\b0\b1\b2&(\r0,\r1,\b1|\r2,\b2|\b0)\cr
\r0\r1\b0\r2\b1\b2&(\r0,\r1,\b1|\b0,\b2|\r2)\cr
\r0\b0\r1\r2\b1\b2&(\r0,\r1,\b1|\b0,\b2|\r2)\cr
\r0\b0\r1\b1\r2\b2&(\r0,\r1,\b1|\b0,\b2|\r2)\cr
\r0\b0\b1\r1\r2\b2&(\r0,\b1,\r1|\b0,\b2|\r2)\cr
\r0\b0\b1\r1\b2\r2&(\r0,\b1,\r1|\b0,\b2|\r2)\cr
\b0\r0\b1\r1\b2\r2&(\b0,\b1,\r1|\r0,\b2|\r2)\cr
\b0\r0\b1\b2\r1\r2&(\b0,\b1,\r1|\r0,\b2|\r2)\cr
\b0\b1\r0\b2\r1\r2&(\b0,\b1,\r1|\r0,\b2|\r2)\cr
\b0\b1\b2\r0\r1\r2&(\b0,\b1,\r1|\b2,\r0|\r2)\cr
}}
$$
Traces of the tableau switching can be already be found find in \ref{Haiman
mixed insertion} and in the \foreign{tables de promotion} of
\ref{Schutzenberger promotion}. The symmetry of our definition with respect to
$i$~and~$j$ has an obvious but important consequence.

\proclaim Theorem. \tabswsym
Tableau switching is involutive: if $X(S,T)=(T',S')$, then $X(T',S')=(S,T)$.
\QED

One particular way to move from the lattice path via $[l,m]$ to the one
via~$[k,n]$, i.e., from the shuffle in which all~$a_i$ precede all~$b_j$ to
the one in which all~$b_j$ precede all~$a_i$, is to start transposing the
final element~$a_{l-1}$ of~$A$ with all elements of~$B$, then to do the same
with~$a_{l-2}$, etc. If, while processing each~$a_i$, we temporarily view the
square labelled by it as ``empty'', then each interchange of labels
corresponds to sliding a label~$b_j$ up or to the left or into the empty
square; the description we so obtain for the transformation of the tableau
labelled by~$B$ before and after processing~$a_i$ is exactly that of in inward
\foreign{jeu de taquin} slide into the square initially labelled by~$a_i$
(see~\ref{Schutzenberger CRGS}, or for instance \ref{Fulton}). It follows
that if $X(S,T)=(T',S')$, then $T'$ is obtained from~$T$ by a sequence of
inward jeu de taquin slides, which we shall write as $T\slid T'$. By symmetry
its follows that one also has $S\dils S'$, i.e., $S'$ is obtained from~$S$ by
a sequence of outward jeu de taquin slides. The following simple observation
will allow us to define tableau switching for skew semistandard tableaux as
well as for skew standard tableaux.

\proclaim Proposition. \wtjdt
If $T\slid T'$, then $T'$ is compatible with a weight~$\alpha$ if and only if
$T$~is.

\proof
Since compatibility of a skew standard tableau with a weight is determined by
the positions of squares being added at successive steps in the chain of
partitions, the general case reduces to the one where $X(S,T)=(T',S')$ with
$T$ of size~$2$ and $S$ of size~$1$. The validity in this case is fairly
obvious by inspection, but to give a formal argument: compatibility
with~$\alpha$ can be tested by comparing the diagonals containing the two
squares of~$T$, and the square of~$S$ cannot be on either of these diagonals;
since the corresponding squares of~$T'$ have at most moved by one diagonal and
certainly cannot have switched positions, the ordering of their diagonals is
unchanged by the slide.
\QED

\proclaim Definition. \sstsdef
Let $T_1,T_2$ be skew semistandard tableaux such that tableau switching can be
applied to the pair $(S_1,S_2)$ of their standardisations; then tableau
switching can also be applied to $(T_1,T_2)$, resulting in a pair of skew
semistandard tableaux $(T'_2,T'_1)=X(T_1,T_2)$ determined by the conditions
that their standardisations are $(S'_2,S'_1)=X(S_1,S_2)$, and the weight
of\/~$T'_i$ is equal to that of\/~$T_i$ ($i=1,2$). The notion of jeu de taquin
slides is also extended to this case; we write $T_2\slid T'_2$ and
$T_1\dils T_1'$.

\heading"Remark"
In fact even more is true: $T'_i$ is $\nu/\kappa$-dominant if and only if
$T_i$~is. This is more difficult to prove than proposition~\wtjdt, although
the proof is still straightforward (one has to prove that after modifying the
companion tableau of~$T_i$ to reflect the slide to~$T'_i$, it still satisfies
the tableau condition). The statement is essentially that of
\ref{pictures, theorem~5.1.1}; see also subsection~3.4 below.

We note one further property of tableau switching that is immediate from its
definition. Given tableaux $T\in\StT{\mu/\nu}$ and~$U\in\StT{\\/\mu}$, a skew
standard tableau of shape~$\\/\nu$ can be formed, which we shall denote
by~$T|U$, by joining together the chains of partitions. In terms or labelled
diagrams, this means we simply combine the labellings of $T$~and~$U$, after
making sure (by adding some offset) that all labels of~$U$ exceed those
of~$T$. Then by similarly joining tableau switching families we get:

\proclaim Proposition. \switchglue
If $X(S,T)=(T',S')$ and $X(S',U)=(U',S'')$, then $X(S,T|U)=(T'|U',S'')$.
\QED

\subsection Jeu de taquin equivalence and dual equivalence.

Jeu de taquin defines an equivalence relation on skew standard tableaux, and
on skew semistandard tableaux, generated by the relations $T\slid T'$; this
relation is called jeu de taquin equivalence. Tableau switching allows us to
define another equivalence relation, called dual equivalence (see~\ref{Haiman
dual equivalence}).

\proclaim Definition. \deqdef
Two skew (semi)standard tableaux $S_1,S_2$ of equal shape are called dual
equivalent if for any tableau~$T$ of appropriate shape, and with
$X(S_i,T)=(T'_i,S'_i)$ for $i=1,2$, one has $T'_1=T'_2$.

It follows from proposition~\switchglue\ that the tableaux $S'_1$~and~$S'_2$
in this definition are also dual equivalent. The meaning of dual equivalence
can be expressed in two ways in terms of jeu de taquin. Firstly, if
$S_1$~and~$S_2$ are dual equivalent, then they have the same shape, and this
remains true whenever the same sequence of outward jeu de taquin slides is
applied to each of them (i.e., each slide starts with the same empty square
in both cases); this is because the tableaux $T'_1$~and~$T'_2$ record the
squares in which the successive slides end. Secondly, when $S_1$~and~$S_2$ are
used to determine sequences of inward slides, then these sequences will always
have the same effect when applied to any tableau~$T$.

It is somewhat surprising that there are any non-trivial instances of dual
equivalence at all. To state the elementary cases, we need some notation for
skew standard tableaux of size~$3$. As mentioned above, a skew standard
tableau~$T\in\StT{\\/\mu}$ can be described by an ordering the squares
of~$\Ylm$; if $|\\/\mu|=3$ those squares lie on distinct diagonals, which are
ordered by index, so we may specify~$T$ using a permutation of the letters
$a$,~$b$,~$c$, which we shall call the type of~$T$. The convention used is
that the three positions in the word correspond to the diagonals, while
$a$,~$b$, and~$c$, respectively indicate the first second and third square
added. As an example, the tableau represented
by~$\smallsquares\Skew(1:1,2|0:3)$ is of type~$cab$.

\proclaim Proposition. \fundamdeq
Let $\\/\mu$ be a shape with $|\\/\mu|=3$. Then $\StT{\\/\mu}$ contains a
tableau of type~$bca$ if and only if it contains a tableau of type~$acb$, and
if so, the two tableaux are dual equivalent. Similarly, $\StT{\\/\mu}$ contains
a tableau of type~$bac$ if and only if it contains one of type~$cab$,
and if so, the two are dual equivalent.

\proof
Note that for tableaux of types in either of the sets $\{bca,acb\}$ or
$\{bac,cab\}$, the relative order among squares of~$\Ylm$ on successive
diagonals is fixed, the only difference being the relative order of the
squares on the first and last diagonal. This observation implies the
statements about the existence of tableaux of these types. The proofs of dual
equivalence will be by induction on the size of the tableau~$T$ in
definition~\deqdef, the case of size~$0$ being trivial. Let
$S_1,S_2\in\StT{\\/\mu}$ and let $T$ be of size~$1$ occupying a square~$t$,
with $X(S_i,T)=(T'_i,S'_i)$ defined for $i=1,2$. Suppose first that~$t$ lies
on a diagonal that does not meet~$\Ylm$; in this case we shall establish the
induction hypothesis (using proposition~\switchglue), by showing that
$T'_1=T'_2$, and that $S'_i$ is of the same type as~$S_i$ for $i=1,2$. Viewing
the computation of $X(S_i,T)$ as three successive applications of an inward
slide to~$T$, this will follow from the fact that the sequence of actual moves
of the label of~$T$ are identical for $i=1,2$ (but possibly occurring during
different slides): the labels of~$S_i$ are restrained to moves to an adjacent
diagonal, and the type cannot change. But the only difference between the
cases $i=1,2$ is the order of the slides into the squares of the outermost
diagonals of~$\Ylm$, and these cannot both be adjacent to that of~$t$, whence
at least one of these slides involves no move (both for $i=1$ and $i=2$); then
the relative order of these two slides makes no difference, and we are done
for this case. The remaining case that $t$ lies on the same diagonal as one of
the squares of~$\Ylm$ can arise only if that is the middle diagonal, and if
$\Ylm$ has the form $\smallsquares\Young(,|)$. Since this is impossible
for~$S_i$ of types $bca$,~$acb$, the proof for that case is therefore
complete. We are left with the case that can be depicted as
$S_1={\smallsquares\Young(1,3|2)}\dils{\smallsquares\Skew(1:1|0:2,3)}=S'_1$
and
$S_2={\smallsquares\Young(1,2|3)}\dils{\smallsquares\Skew(1:2|0:1,3)}=S'_2$;
we see that $T'_1=T'_2$, and that $S'_1$ and~$S'_2$ are of types
$bca$~and~$acb$, whence they are dual equivalent by the case just completed.
\QED

This basic case implies many others by the following immediate consequence
proposition~\switchglue:

\proclaim Proposition. \gluedeq
If $S_1$~and~$S_2$ are dual equivalent, and also $T_1$~and~$T_2$, then so
are $S_1|T_1$~and~$S_2|T_2$.
\QED

\subsection Confluence of jeu de taquin.

A crucial property of jeu de taquin is its confluence, i.e., the property that
whenever $S\slid T_1$ and $S\slid T_2$, then there exists a tableau~$U$ such
that $T_1\slid U$ and~$T_2\slid U$. Since any sequence of slides can be
extended until a tableau of partition shape is reached, and no further (Young
tableaux are the normal forms for jeu de taquin), this is equivalent to saying
that for any skew tableau~$S$ there is a \emph{unique} Young tableau~$P$ such
that $S\slid P$. Although this property will follow as a corollary to our main
theorem in~\Sec3, it is so important that it is worth while to give an
independent proof here, that does not require any further constructions. The
proof that we shall present was inspired by the proof given in~\ref{Kimmo's
thesis} (although it must be admitted that by translating the reasoning using
tableau switching, it has become rather similar to the arguments contained
in~\ref{Haiman dual equivalence}); we have in fact already established the
essential part of the argument as proposition~\fundamdeq. By our second
description of dual equivalence above, confluence of jeu de taquin amounts to
the following.

\proclaim Theorem. \Ytdeq
All standard Young tableaux of shape~$\\$ are dual equivalent.

\proof
The proof is by induction on~$|\\|$; the case $|\\|=0$ is trivial (in fact all
cases with $|\\|\leq3$ are either trivial or already established). By the
induction hypothesis and proposition~\gluedeq, one has dual equivalence among
all members of any subclass~$C(\\,s)$ of all standard Young tableaux of
shape~$\\$ for which the highest entry occupies a given square~$s$ (as chains
of partitions, they share a predecessor of~$\\$ in~$\Ylat$). It will suffice
to establish dual equivalence between some pair of representatives of any two
such subclasses $C(\\,p)$, $C(\\,r)$; we shall assume that the diagonal of~$p$
has smaller index than that of~$r$. The squares $p$~and~$r$ are both corners
of the Young diagram~$Y(\\)$, and the Young diagram $Y(\\)\setminus\{p,r\}$
contains at least one corner~$q$ that lies on a diagonal with index between
those of $p$~and~$q$ (e.g., its corner in the row above~$p$). Let
$\mu\in\Part$ be such that $Y(\mu)=Y(\\)\setminus\{p,q,r\}$, and let
$S_1,S_2\in\StT{\\/\mu}$ be the tableaux of respective types $cab$~and~$bac$;
these tableaux are dual equivalent by proposition~\fundamdeq. For any
$R\in\StT\mu$ one has $R|S_1\in C(\\,p)$ and $R|S_2\in C(\\,r)$, and by
proposition~\gluedeq, $R|S_1$~and~$R|S_2$ are dual equivalent.
\QED

\heading "Remark"
The statement of the theorem corresponds to a global form of confluence of jeu
de taquin: any two sequences of jeu de taquin slides starting with the same
skew tableau~$S$ will, when extended sufficiently far, always end in the same
normal form (Young tableau). We have however derived this from a local form of
confluence: if two different jeu de taquin slides are applied to~$S$, then the
resulting tableaux can be made equal by applying at most two more slides to
each of them. The fact that this local confluence is achieved with the same
number of slides along either route avoids difficulties in deriving global
confluence that otherwise exist (see~\ref{Kimmo's thesis}); it is therefore
essential that our definition of jeu de taquin ensures, by including certain
slides that change only the skew shape of a tableau but not the associated
diagram, that all sliding sequences between two given tableaux
always have the same length.

By definition, two tableaux are dual equivalent if their equality of shape is
preserved under sequences of outward slides, which means also that the two
sequences of inward slides (applied to some tableau~$T$) determined by them
always have the same effect. One may ask whether the same properties also hold
with ``inward'' and ``outward'' interchanged (indeed this is a requirement for
dual equivalence in the original definition in~\ref{Haiman dual equivalence}).
To prove that this is the case, we shall use the fact that the sub-poset of
the Young lattice of all partitions contained in some fixed rectangular
partition~$\rho$, possesses an involutive poset anti-isomorphism. If $\rho$ is
the partition with $n$ non-zero parts, all equal to~$m$, then this
anti-isomorphism, which we shall denote by~$\\\mapsto\\^\diamond$ for
$\\\subset\rho$, satisfies $\\^\diamond_i=m-\\_{n-1-i}$ for $i\in\upto n$. For
$T\in\StT{\\/\mu}$ with $\\\subset\rho$, we also define
$T^\diamond\in\StT{\mu^\diamond/\\^\diamond}$ by applying the anti-isomorphism
to all partitions in the chain of~$T$, and reversing their order. By a similar
operation applied to tableau switching families, it can be seen that
$X(S,T)=(T',S')$ implies
$X(T^\diamond,S^\diamond)=({S'}^\diamond,{T'}^\diamond)$.

\proclaim Proposition. \revdeq
Two tableaux $S_1,S_2$ of the same shape are dual equivalent if and only if
for any tableau~$T$ of appropriate shape, and with $X(T,S_i)=(S'_i,T'_i)$ for
$i=1,2$, one has $T'_1=T'_2$.

\proof
Let us temporarily call the relation in the second half of the proposition
reverse dual equivalence. Suppose first that $S_1$~and~$S_2$ are reverse dual
equivalent and of shape~$\\/\mu$. Let $T\in\ST(\mu/0)$, and put
$(S'_i,T')=X(T,S_i)$ for $i=1,2$; then $S'_1$~and~$S'_2$ are dual equivalent
by theorem~\Ytdeq, so $S_1$~and~$S_2$ are also dual equivalent by the remark
after definition~\deqdef. (In fact theorem~\Ytdeq\ is equivalent to the ``if''
part of the current proposition, since Young tableaux of equal shape are
trivially reverse dual equivalent.) Conversely, suppose that $S_1$~and~$S_2$
are dual equivalent; then for any choice of a rectangular
partition~$\rho\supset\\$, the tableaux $S_1^\diamond$~and~$S_2^\diamond$ are
reverse dual equivalent, whence (by what we just proved) they are dual
equivalent, and so $S_1=S_1^{\diamond\diamond}$ and
$S_2=S_2^{\diamond\diamond}$ are reverse dual equivalent.
\QED

This proposition provides us with an effective test of dual equivalence. Given
two tableaux $S_1,S_2$ of equal shape, one successively applies to both
tableaux inward slides into the same squares; if at any moment their shapes
become different, then $S_1$~and~$S_2$ are not (reverse) dual equivalent, but
if the tableaux become Young tableaux without exhibiting shape difference,
then one has established $X(S'_i,T')=(T,S_i)$ for $i=1,2$ with
$S'_1$~and~$S'_2$ dual equivalent (Young) tableaux, so $S_1$~and~$S_2$ are
dual equivalent. For jeu de taquin equivalence one has of course also a test,
thanks to confluence, which consists of independently reducing both tableaux
to Young tableaux, which will be equal in case of jeu de taquin equivalence.

\proclaim Corollary. \dubleqcor
If $S_1,S_2$ are both jeu de taquin equivalent and dual equivalent, then
$S_1=S_2$.
\QED

\subsection An alternative Littlewood-Richardson rule.

Let $\chi\in\Sk$; the fact that to each skew tableau $T\in\StT\chi$ there is
associated a unique Young tableau~$P$ obtainable from it by jeu de taquin,
allows us to subdivide $\StT\chi$ according to the shape of this~$P$. For
$\nu\in\Part_{|\chi|}$ define
$\StT\chi\snu=\Union_{P\in\StT\nu}\StT\chi^{\slid P}$, where $\StT\chi^{\slid
P}=\setof T\in\StT\chi:T\slid P\endset$ for any $P\in\StT\nu$. All fibres
$\StT\chi^{\slid P}$ of the map $\StT\chi\snu\to\StT\nu$, sending $T\mapsto P$
when $T\slid P$ have the same number of elements; in fact there are canonical
bijections between the fibres, so that $\StT\chi\snu$ is in natural bijection
with the Cartesian product of~$\StT\nu$ and any one of the fibres:

\proclaim Proposition. \LRjdtST
Let $\\/\mu\in\Sk$ and $\nu\in\Part$, and let $C\in\StT\nu$ be fixed. Then
there is a bijection $\phi\:\StT{\\/\mu}^{\slid C}\times\StT\nu\to
\StT{\\/\mu}\snu$, that can be characterised by the conditions that
$\phi(L,P)\slid P$ and that $\phi(L,P)$ is dual equivalent to~$L$.

\proof
We can construct $\phi(L,P)$, which is uniquely determined due to
corollary~\dubleqcor, as follows. Choose any $Q\in\StT\mu$, then
$X(Q,L)=(C,S)$ for some $S\in\StT{\\/\nu}$; moreover $X(P,S)=(Q,T)$ for some
$T\in\StT\chi^{\slid P}$ by theorem~\Ytdeq, since $S\slid Q$, and we set
$\phi(L,P)=T$. It is clear that $\phi(L,P)\slid P$, while $L$ and
$T=\phi(L,P)$ are dual equivalent since the Young tableaux $C$~and~$P$ are.
The bijectivity of~$\phi$ follows from the fact that $P$, $S$, and~$L$ are
readily reconstructed from~$T$.
\QED

As an example of this construction, let
$$ \def~{\mkern-2mu }
   C={\smallsquares\Young(0,1,2,7,1~1|3,5,9|4,6,1~0|8,1~2)}, \qquad
   L={\smallsquares\Skew(4:2,7,1~1|3:5,9|1:0,1,1~0|1:3,1~2|0:4,6|0:8)}, \qquad
   P={\smallsquares\Young(0,2,4,5,9|1,3,7|6,1~0,1~2|8,1~1)};
$$
choosing $Q={\smallsquares\Young(0,1,5,8|2,4,7|3|6)}$
one gets $S={\smallsquares\Skew(5:1,8|3:2,5|3:7|2:4|0:0,6|0:3)}$,
and $X(P,S)=(Q,T)$, where
$ \def~{\mkern-2mu }
   T=\phi(L,P)=
   {\smallsquares\Skew(4:2,5,9|3:4,7|1:0,3,1~0|1:1,1~2|0:6,1~1|0:8)}
$.
\smallskip\noindent
The reader may check that the same tableau $\phi(L,P)$ is obtained for other
choices of~$Q$. The above proposition can be generalised to skew semistandard
tableaux, using proposition~\wtjdt. It suffices to replace $\StT\nu$
and~$\StT{\\/\mu}$ respectively by $\Tab(\nu;n)$ and $\Tab(\chi;n)$, and
$\StT{\\/\mu}\snu$ by $\Tab(\chi;n)\snu$, which is defined as
$\Union_{C\in\Tab(\nu;n)}\JdTt\chi,C$ where $\JdTt\chi,C=\setof
L\in\Tab(\chi;n):L\slid C\endset$. We obtain the following corollary, either
by applying proposition~\LRjdtST, or by reusing its proof almost literally.

\proclaim Corollary. \LRjdtcor
Let $\chi\in\Sk$, $\nu\in\Part$, $n\in\N$, and fix $C\in\Tab(\nu;n)$. Then
there is a bijection
$\phi\:\JdTt\chi,C\times\Tab(\nu;n)\to\Tab(\chi;n)\snu$, that can be
characterised by the conditions that $\phi(L,P)\slid P$ and that $\phi(L,P)$
is dual equivalent to~$L$.
\QED

If one combines the inverses of such bijections~$\phi$ for all
$\nu\in\Part_{|\chi|,n}$, one obtains a bijection~$\Rob'$ defined on
$\Union_{\nu\in\Part_{|\chi|,n}}\Tab(\chi;n)\snu=\Tab(\chi;n)$, whose codomain
strongly resembles that specified for Robinson's bijection~$\Rob$ at the right
hand side of~(\Robinson); the difference is that, in each component of the
union over $\nu\in\Part_{|\chi|,n}$, the factor $\LR(\chi,\nu)$ is replaced by
$\JdTt\chi,C$, for some~$C\in\Tab(\nu;n)$ chosen separately for every~$\nu$.
Writing $\Rob'(T)=(\Rob'_0(T),\Rob'_1(T))$, the map $\Rob'_1$ preserves weight
(since $T\slid\Rob'_1(T)$); therefore the bijection~$\Rob'$ gives us an
alternative expression for the value of Littlewood-Richardson coefficients:

\proclaim Corollary. \altLRR
If the skew shape $\chi$ represents~$\\*\mu$, then the Littlewood-Richardson
coefficient~$c_{\\,\mu}^\nu$ equals $\Card{\JdTt\chi,C}$ for
any~$C\in\Tab(\nu)$. More generally, $\Card{\JdTt\chi,C}=c_\chi^\nu$
for any skew shape~$\chi$.
\QED

Although this version of the Littlewood-Richardson rule has the advantage of
already being proved, it is of little practical use in its present form, since
the sets $\JdTt\chi,C$ cannot be enumerated in any useful manner. However, we
shall in the next section prove the identity
$$ 
  \JdTt\chi,{\1_\nu}=\LRt\chi,\nu; \label(\LRfibereq)
$$
it follows that for the special choice $C=\1_\nu\in\Tab(\nu;n)$ for all~$\nu$,
the bijection~$\Rob'$ exactly matches the specification of Robinson's
bijection, providing a proof of the Littlewood-Richardson rule. Indeed, we
shall define $\Rob$ so that it coincides with this specialisation of~$\Rob'$.
The remark we made following definition~\sstsdef\ implies~(\LRfibereq), so we
could complete our proof of the Littlewood-Richardson rule by proving that
remark. Using coplactic operations, as we shall do, is certainly not the
simplest way to prove~(\LRfibereq), but it provides a better insight into
Robinson's bijection, and in particular it gives a simpler description
of~$\Rob_0$ than the one that can be extracted from what has been presented so
far. It will also lead to a proof of the Littlewood-Richardson rule that does
not depend on any of the non-trivial results derived in this section.

By proposition~\LRswitchprop, (\LRfibereq) also implies
$\Card{\JdTt\chi,{\1_\nu}}=\Card{\JdTt\nu/\mu,{\1_\\}}$ when $\chi$
represents~$\\*\mu$, and hence $c_{\\,\mu}^\nu=c_{\nu/\mu}^\\$. On the other
hand $c_{\nu/\mu}^\\=c_{\nu/\\}^\mu$ can be obtained without
using~(\LRfibereq): the relation $X(\1_\mu,T)=(\1_\\,T^*)$ defines a~bijection
between tableaux $T\in\JdTt\nu/\mu,{\1_\\}$ and
$T^*\in\JdTt\nu/\\,{\1_\mu}$ (this can be generalised by replacing
$\1_\mu$~and~$\1_\\$ by other fixed tableaux of the same shape). For instance,
here is the computation for the tableau~$L^*$ corresponding to the
Littlewood-Richardson tableau~$L$ of~(\LRtablex).
$$\def\r#1{{\it#1\/}}\def\b#1{{\bf#1}}
L= {\maybesmall\Skew(4:0,0|2:0,1,1|1:0,1,2,2|0:0,1,3|0:1,4)}
,\quad
 {\maybesmall
  \Young(\r0,\r0,\r0,\r0,\b0,\b0
	|\r1,\r1,\b0,\b1,\b1
	|\r2,\b0,\b1,\b2,\b2
	|\b0,\b1,\b3
	|\b2,\b4)
 }
\buildrel X \over \longleftrightarrow
 {\maybesmall
  \Young(\b0,\b0,\b0,\b0,\b0,\r0
	|\b1,\b1,\b1,\b1,\r0
	|\b2,\b2,\b2,\r0,\r1
	|\b3,\r0,\r2
	|\b4,\r1)
 }
,\quad
L^*={\maybesmall\Skew(5:0|4:0|3:0,1|1:0,2|1:1)}
\nn
$$
The identity $c_{\nu/\mu}^\\=c_{\nu/\\}^\mu$ is of course related to the
symmetry $c_{\\,\mu}^\nu=c_{\mu,\\}^\nu$ that is obvious from the definition.
Writing $\LRt\\*\mu,\nu$ for $\LRt\chi,\nu$ when $\chi$ represents $\\*\mu$,
the following description of a bijection between $\LRt\\*\mu,\nu$
and~$\LRt\mu*\\,\nu$ can be found, assuming~(\LRfibereq): given
$L\in\LRt\\*\mu,\nu$, determine $\bar L_\\\in\LRt\nu/\mu,\\$ corresponding to
it by proposition~\LRswitchprop\ (i.e., the companion tableau of the
subtableau~$L_\\$ of~$L$), then determine $\bar L_\\^*\in\LRt\nu/\\,\mu$ as
above by $X(\1_\mu,\bar L_\\)=(\1_\\,\bar L_\\^*)$, and finally apply the
bijection corresponding to proposition~\LRswitchprop\ in the opposite
direction to~$\bar L_\\^*$ so as to obtain a tableau in~$\LRt\mu*\\,\nu$. It
does not appear that this somewhat complicated process can be simplified.

\endofchapter

 \par

\ifx\macrosloaded\undefined
  \input marxmax
  

\newif\ifnarrow
\ifdim \hsize<14.5cm \narrowtrue \else\narrowfalse \fi

\def\maybesmall{\ifnarrow\smallsquares\fi}

\def\ol#1
{\vbox{\vskip0.3pt \hrule \vskip 0.8pt \hbox{$#1$}}}

\def\ul#1
{\vtop{\hbox{$#1$}\vskip 0.8pt \hrule \vskip0.3pt }}

\font\frak=eufm10
\font\csc=cmcsc10
\font\URLfont=cmtex10 
\hyphenchar\URLfont=-1 

\mathorddef\Part=P2 
\mathorddef\A=A0 \mathorddef\B=B0 
\def\GL{\boldmath{GL}}
\def\Ylat{(\Part,{\subset})}
\def\upto#1{[#1]}
\opdef\Pic
\altopdef\Tab{SST}
\opdef\ST
\opdef\LR
\opdef\JdT
\opdef\wt
\mathorddef\T=T2 
\mathorddef\Sk=S2 
\def\SFn{\Lambda_n}
\def\Can#1{\boldmath1_{#1}}
\def\Ylm{Y(\\/\mu)}
\def\bT{\bar T}
\def\LRt#1,#2{\LR(#1,#2)}
\def\JdTt#1,#2{\Tab(#1)^{\slid{#2}}}
\mathorddef\Rob=R2 
\def\StT#1{\ST(#1)}
\def\snu{^{\slid\nu}}
\def\Uqgln{U_q(\hbox{\frak gl}_n)}
\def\leqr{\leq_r}
\def\leqKan{\leq_{\romath K}}

\def\wKan{w_{\romath K}}
\def\wSem{w_{\romath S}}

\mathdef\slid=Rel \triangleright
\mathdef\dils=Rel \triangleleft

\long\def\quotation#1
{{\smallskip\narrower\it\noindent\llap{``\/}#1\unskip''\smallskip}\noindent}

{\catcode`\.=\active\catcode`\/=\active\catcode`\_=\active
 \gdef\URL
 {\begingroup\def~{\char`\~}%
  \catcode`\.=\active\catcode`\/=\active\catcode`\_=\active
  \let_=\_%
  \def.{\char`\.\penalty\hyphenpenalty}\def/{\char`\/\penalty\hyphenpenalty}%
  \URLfont\URLnext
}}
\def\URLnext#1{http://#1\endgroup}

\hyphenation
{semi-standard Yamano-uchi mono-mial mono-mials inter-related combi-na-torics
 combi-na-torial algo-rithm valid-ity}

  \writelab{}
  \def\jobname{lrr}
\readrefs
  \let\endofchapter=\bye
  \input lrr1.lab
  \input lrr2.lab
\fi

\newsection Coplactic operations.

In this section we shall introduce another kind of operations on skew
semistandard tableaux, which we shall call \newterm{coplactic operations}.
Unlike jeu de taquin, these transformations do not change the shape of a
tableau, but rather its weight, by changing the value of one of the entries.
The basic definitions can be formulated most easily in terms of finite words
over the alphabet~$\upto n$. In the application to tableaux, these words will
be the ones obtained by listing the entries of skew semistandard tableaux
using a valid reading order as described in~\Sec\readingsec; this means in
particular that notions of ``left'' and ``right'' within words will
(unfortunately) get a more or less opposite interpretation within tableaux.

\subsection Coplactic operations on words.

We first fix some terminology pertaining to words. A \newterm{word} over a
set~$A$ (called the alphabet) is a finite (possibly empty) sequence of
elements of~$A$, arranged from left to right; the elements of the sequence
forming a word~$w$ are called the \newterm{letters} of~$w$. The set of all
words of length~$l$ over~$A$ is denoted by~$A^l$, and
$A^*=\Union_{l\in\N}A^l$. The concatenation of two words $u,v\in A^*$ will be
denoted by~$uv$ (the sequence~$u$ of letters, followed by the sequence~$v$);
this defines an associative product on~$A^*$. Whenever a word~$w$ can be
written as~$uv$, the word~$u$ is called a \newterm{prefix} of~$w$, and $v$ a
\newterm{suffix} of~$w$; a \newterm{subword} of~$w$ is any word that can be
obtained by removing from~$w$ a (possibly empty) prefix and a suffix.

A word $w\in\upto n^*$ will be called \newterm{dominant} for~$i\in\upto{n-1}$
if every prefix of~$w$ contains at least as many letters~$i$ as letters~$i+1$;
similarly $w$ will be called \newterm{anti-dominant} for~$i$ if every suffix
of~$w$ contains at least as many letters~$i+1$ as letters~$i$. If~$w$ is both
dominant and anti-dominant for~$i$, it will be called \newterm{neutral}
for~$i$. It is easy to see that when $w$ is either dominant or anti-dominant
for~$i$, then it is neutral for~$i$ if and only if it contains exactly as many
letters~$i$ as letters~$i+1$. In a word~$w$ that is neutral for~$i$ there is a
matching between letters~$i$ and letters~$i+1$, like properly matched left and
right parentheses (the words in $\{i,i+1\}^*$ that are neutral for~$i$
constitute a so-called Dyck language), and for such~$w$, any word of the form
$uwv$ will be dominant (respectively anti-dominant, or neutral) for~$i$ if and
only if $uv$ is so.

\proclaim Definition. \coplactdef
A coplactic operation in $\upto n^*$ is a transition between words $w=uiv$
and $w'=u(i+1)v$ where $u,v\in\upto n^*$ and $i\in\upto{n-1}$ are such that
$u$ is anti-dominant for~$i$ and $v$ is dominant for~$i$. In this case we
shall write $w=e_i(w')$ and $w'=f_i(w)$.

For instance, in the following word over the alphabet~$\upto6$, decrementing
by~$1$ any one of the numbers with an underline, or incrementing by~$1$ any one of
the numbers with an overline, constitutes a coplactic operation, and no other
coplactic operations are possible; the word is dominant for~$1$ and neutral
for~$2$.
$$ \let\e\underline \let\f\overline \let~\enspace
  \e4~0~1~5~2~\e{\f1}~3~\e5~0~1~\f4~2~\f0~0~1~2~\f3~3~4
  \label(\copex)
$$
In the definition, the letter~$i$ in~$w$ that is changed to~$i+1$ in~$w'$ is
not contained in any subword~$u_0iv_0$ of~$w$ that is anti-dominant for~$i$,
for $v_0$ would then have strictly more letters~$i+1$ than letters~$i$,
contradicting the dominance for~$i$ of~$v$; similarly the indicated
letter~$i+1$ in~$w'$ is not contained in any subword that is dominant for~$i$.
In particular the changing letters are not contained in any subword that is
neutral for~$i$, and $u$ is the longest prefix of~$w$ that is anti-dominant
for~$i$, while $v$ is the longest suffix of~$w'$ that is dominant for~$i$.
Hence the expressions $f_i(w)$ and $e_i(w')$, when defined, are unambiguous.

\proclaim Proposition. \copdefdprop
The expression $e_i(w)$ is defined unless $w$ is dominant for~$i$, and
$f_i(w)$ is defined unless $w$ is anti-dominant for~$i$.

\proof
We shall proof the latter statement, the proof of former being analogous. Let
$u$ be the longest prefix of~$w$ that is anti-dominant for~$i$; clearly
$f_i(w)$ cannot be defined if $u=w$. Otherwise $w=uiv$ for some~$v$, and we
show by induction on its length that~$v$ is dominant for~$i$, which will prove
the proposition. The cases where $v$ is empty or ends with a letter other
than~$i+1$ are trivial, so assume $v=v'(i+1)$ and suppose $v$ is not dominant
for~$i$ while (by induction) $v'$~is. Then $v'$ has as many letters~$i$ as
letters~$i+1$, and is therefore neutral for~$i$, so that $w=uiv'(i+1)$ is
anti-dominant for~$i$ since $ui(i+1)$ is; a contradiction.
\QED

We associate a weight $\wt w\in\N^n$ to words $w\in\upto n^*$ in the same way
as to semistandard tableaux, i.e., the component~$i$ of~$\wt w$ counts the
number of of letters~$i$ in~$w$. Then when $w=e_i(w')$, one
has~$\wt w\succ\wt w'$; for this reason the operations~$e_i$ are called
\newterm{raising operations}, and the~$f_i$ are called \newterm{lowering
operations}. Starting with any $w\in\upto n^*$ one can iterate application of
a fixed~$e_i$ until, after a finite number of iterations, $w$ is transformed
into a word that is dominant for~$i$. More generally any sequence of
applications of operations~$e_i$, where $i$ is allowed to vary, must
eventually terminate, producing a word that is simultaneously dominant for
all~$i\in\upto{n-1}$. Such a word will be simply called \newterm{dominant}, or
in older terminology a \newterm{lattice permutation}, which means that the
weight of any prefix is a partition (of the length of the prefix into
$n$~parts). Observe the similarity with $0$-dominance for semistandard
tableaux, as characterised in proposition~\LRlatticechar; this is what
motivated our choice of terminology. For instance, for the word in~(\copex),
if we choose at each step to apply the operation~$e_i$ with minimal
possible~$i$, the sequence of operations applied is
$\let~=\>{e_0,~e_3,e_2,e_1,e_0,~e_4,e_3,e_2,e_1,~e_4,e_3,e_2,e_1}$, which
respectively decrease the letters at positions
$\let~=\>{5,~0,0,0,0,~7,10,10,11,~3,3,3,4}$ from the left, leading finally
to the dominant word $\let~=\,0~0~1~2~1~0~3~4~0~1~2~1~0~0~1~2~3~3~4$. Thus the
raising operations~$e_i$ define a rewrite system on~$\upto n^*$, whose normal
forms are the dominant words. We shall show below that this rewrite system is
confluent; for instance the dominant word just obtained from the word
in~(\copex) can also be obtained by always applying the~$e_i$ with
maximal possible~$i$, which leads to the sequence
$\let~=\>{e_4,~e_3,~e_3,e_4,~e_2,~e_2,e_3,~e_1,e_2,~e_0,e_1,~e_0,e_1}$,
respectively affecting letters at positions
$\let~=\>{7,~10,~0,3,~10,~0,3,~0,3,~5,11,~0,4}$.

The coplactic operations define a labelled directed graph on the set $\upto
n^*$, with an edge labelled~$i$ going from $w$~to~$w'$ whenever $f_i(w)=w'$;
we shall call this the \newterm{coplactic graph} on~$\upto n^*$. For each
$w\in\upto n^*$, we shall call the connected component of the coplactic graph
on~$\upto n^*$ containing~$w$ the coplactic graph of~$w$; we consider this to
be a rooted graph, with as root the element~$w$ itself. For $n=2$, these
coplactic graphs are linear with a dominant word at one end and an
anti-dominant word at the other end; distinct raising operations do not
commute however, so that for $n>2$ the coplactic graph associated to~$w$ can
contain other words with the same weight (for instance the coplactic graph
of~$1\,0\,2$ also contains~$2\,0\,1$). Coplactic graphs are isomorphic to the
crystal graphs for irreducible integrable $\Uqgln$-modules of~\ref{crystal
graphs}; this places their properties in a broader perspective. Their
structure is intricate, and not easy to describe in an independent way;
however, as we shall see, it occurs frequently that distinct words have
isomorphic coplactic graphs (for instance the coplactic graphs of $1\,0\,2$
and~$0\,2\,1$ are isomorphic). The coplactic graph of the word~$w$ in~(\copex)
for $n=6$ contains $53460$~words, $120$~of which have the same weight as~$w$
(for instance $\let~=\,2~0~1~5~3~1~4~5~0~1~3~2~0~0~1~2~4~3~4$); we shall not
attempt to depict this graph. However, we encourage the reader to draw some
small coplactic graphs, e.g., for the word~$0\,1\,0\,1$ and $n=4$ (cf.\
\ref{quantum straightening, Figure~2}).

\subsection Coplactic operations on tableaux.

As we mentioned, our interest in coplactic operations is in applying them to
tableaux rather than to words. For the purpose of defining coplactic
operations, the entries of a tableau will be considered as letters of a word
that have been mapped in some order onto the squares of a skew diagram. In
view of the similarity between dominance of words and $0$-dominance of
tableaux, it should come as no surprise that the order in which the entries of
a tableau are strung together into a word is a valid reading order as
discussed in~\Sec\readingsec. Nevertheless the following properties are quite
remarkable: despite the loss of information about rows and columns of the skew
diagram caused by this reading process, the tableau conditions are preserved
when coplactic operations are applied to the word, regardless of which valid
reading order is used, and in fact the changes to the entries constituting the
coplactic operations are themselves independent of that order, i.e., the same
entry is affected, at whatever place in the word the reading order places it.
For instance, here is a tableau with the coplactic operations that can be
applied to it, labelled as in~(\copex), and two of its similarly labelled
reading words (for the Semitic and Kanji reading order, respectively) that
could be used to determine those possible coplactic operations.
$$
  \let\e\ul \let\f\ol \let~\,
  \Skew(3:0,\f1|1:\f0,\e1,1,3|0:0,2,\f2,\e3|0:1,4,\f4,\e5|0:\f3,5)
\qquad \let\e\underline \let\f\overline 
  \f1~0~3~1~\e1~\f0~\e3~\f2~2~0~\e5~\f4~4~1~5~\f3
\qquad
  \f1~3~0~1~\e3~\e5~\e1~\f2~\f4~\f0~2~4~5~0~1~\f3
$$

We shall now give a more formal definition. A valid reading order
for~$\chi\in\Sk$ is a total ordering~`$\leqr$' on~$Y(\chi)$, such that
$(i,j)\leqr(i',j')$ whenever $i\leq i'$ and~$j\geq j'$. A corresponding map
$w_r\:\Tab(\chi,n)\to\upto n^*$ is defined by $w_r(T)=\fold{}(Ts_0..s_k)$,
where $Y(\chi)=\{\li(s0..k)\}$ with $\fold{<_r}(s0..k)$; in other words, $w_r$
forms a list of all entries of~$T$, in increasing order for~`$\leqr$' of their
squares. Now if some coplactic operation~$c$ can be applied to~$w_r(T)$,
changing a letter~$l$, then the square~$s$ of the entry of~$T$ corresponding
to~$l$ will be called the variable square for this operation~$c$; formally, if
$k$~letters of~$w_r(T)$ precede~$l$, then $\Card{\setof t\in Y(\chi):
t<_rs\endset}=k$. Then if changing the entry of~$s$ in~$T$ in the same way
as~$l$ results in a tableau~$U\in\Tab(\chi,n)$, we define $U=c(T,\leqr)$.
Consequently $w_r(c(T,\leqr))=c(w_r(T))$, when defined.

\proclaim Proposition. \coptabprop
Let $\chi\in\Sk$, let $T\in\Tab(\chi,n)$ and let $c$ be a coplactic operation
$e_i$~or~$f_i$ with $i\in\upto{n-1}$. Then for any valid reading
order~`$\leqr$', the tableau~$c(T,\leqr)$ is defined if and only if
$c(w_r(T))$~is; moreover, this condition and (when it holds) the value
of~$c(T,\leqr)$ do not depend on~`$\leqr$'.

\proof
We shall use a process of reduction: $T\in\Tab(\chi,n)$ is simplified by
successively removing certain sets of squares from $Y(\chi)$, restricting $T$
and~`$\leqr$' to the remainder. At each step the change to $w_r(T)$ will be
the removal of a subword that is neutral for~$i$, so that neither the
condition whether or not $c(w_r(T))$ is defined, nor the variable square are
affected. The reduction steps are of two types. The first type is the removal
of squares whose entries are not $i$~or~$i+1$, as the corresponding subwords
of~$w_r(T)$ (of length~$1$) are neutral for~$i$; this reduces us to the case
that $T$ has entries $i$~and~$i+1$ only. In that case the second type of
reduction applies, which removes a maximal rectangle within~$Y(\chi)$
consisting of two rows of squares (since the columns of~$Y(\chi)$ now have at
most $2$~squares, maximality means the rectangle cannot be extended into the
columns to its left or right). One easily checks that the letters of~$w_r(T)$
corresponding to such a rectangle form a subword that is neutral for~$i$. When
no further reduction of this type is possible, $T$ is reduced to a tableau
with at most one square in any column, and regardless of the original reading
order, $w_r(T)$ lists their entries from the rightmost column to the leftmost
one. Assuming now that $c(w_r(T))$ is defined, it follows from the definition
of coplactic operations for words that the variable square does not have a
left neighbour with entry~$i+1$, nor a right neighbour with entry~$i$; this
shows in the reduced case that changing the entry of the variable square does
not violate weak increase along rows. Neither is it possible that such a
neighbouring entry of the variable square was removed by the second type of
reduction (the variable square itself would then have to have been in a column
of length~$2$, which it was not), showing that weak increase along rows is
preserved in the unreduced case as well; strict increase down columns is
preserved since their are no other entries $i$~or~$i+1$ in the column of the
variable square.
\QED

\proclaim Definition.
On $\Tab(\chi,n)$ the coplactic operations $e_i$ and~$f_i$ (for $i<n$) are
(partially) defined by $e_i(T)=e_i(T,\leqr)$ and $f_i(T)=f_i(T,\leqr)$ for an
arbitrary valid reading order~`$\leqr$'; this is taken to mean also that the
left hand sides are undefined when the corresponding right hand sides are.

We call $T\in\Tab(\chi)$ dominant if $w_r(T)$ is dominant for any (and hence
every) valid reading order~$r$; it means $e_i(T)$ is undefined for
all~$i\in\N$. From proposition~\LRlatticechar\ we get the following
characterisation.

\proclaim Corollary. \LRdomincorr
For any $\chi\in\Sk$, the subset of dominant elements of $\Tab(\chi)$ is equal
to~$\LR(\chi)$.
\QED

\subsection The main theorem: commutation.

As a result of proposition~\coptabprop, one obtains many instances of distinct
words with isomorphic coplactic graphs, namely words obtained as $w_r(T)$ for
fixed~$T$ and different reading orders~`$\leqr$'. In fact we shall presently
find even more instances than can be found in this manner. The situation
resembles to that of jeu de taquin, where we found remarkably many cases of
dual equivalent skew semistandard tableaux. We shall now give a theorem that
explains both these phenomena, and at the same time essentially proves the
Littlewood-Richardson rule. This theorem, which is the only one in our paper
with a somewhat technical proof, simply states that jeu de taquin commutes
with coplactic operations. This implies that coplactic operations transform
tableaux into dual equivalent ones, and that words obtained from jeu de taquin
equivalent tableaux using any reading order always have isomorphic coplactic
graphs.

\proclaim Theorem. \copcommthm
Coplactic operations $e_i$ and~$f_i$ on tableaux commute with jeu de taquin
slides in the following sense. If $S\in\Tab(\mu/\nu)$, \ $T\in\Tab(\\/\mu)$,
and $X(S,T)=(T',S')$, then $e_i(T)$ is defined if and only if $e_i(T')$~is,
and if so, one has $X(S,e_i(T))=(e_i(T'),S')$; the same holds when $e_i$ is
replaced by~$f_i$.

In the course of the proof we shall need to draw some specific configurations
that may occur within the tableaux. These will involve entries $i$~and~$i+1$
only; in order to fit them into the squares, we subtract~$i$ from each,
representing them respectively as $0$~and~$1$ (one might also say the drawings
assume $i=0$).

\proof
Since the standardisation of~$T$ can be written as $U|V|W$ where $U$, $V$, $W$
are the standardisations of the subtableaux of~$T$ of entries less than~$i$,
in $\{i,i+1\}$, and greater than~$i+1$, respectively, we can reduce by
proposition~\switchglue\ to the case that $T$ has entries $i$~and~$i+1$ only;
we may also assume~$|\mu/\nu|=1$. Since $T_0=e_i(T_1)$ means the same as
$f_i(T_0)=T_1$, it suffices to consider the operations~$e_i$. We may use any
valid reading order to determine coplactic operations on tableaux; it will be
convenient to use the Kanji reading order, which we shall denote by~$\leqKan$.
The only differences between $\wKan(T)$ and~$\wKan(T')$ are due to transitions
$\smallsquares\Young(,0|1,1)\to\Young(0,|1,1)$ and
$\smallsquares\Young(0,0|,1)\to\Young(0,0|1,)$ during the slide~$T\slid{T'}$;
the effect of these transitions on~$\wKan(T)$ amounts to interchanging one
letter with a word $i(i+1)$. Since that word is neutral for~$i$, this does not
affect the dominance or anti-dominance for~$i$ of any subword containing the
change. In particular $\wKan(T')$ is dominant for~$i$ if and only
if~$\wKan(T)$~is, whence $e_i(T)$ is defined if and only if $e_i(T')$~is; we
assume henceforth that this is the case. Let~$v$ be the variable square for
the application of~$e_i$ to~$T$, and $v'$ the variable square for the
application of~$e_i$ to~$T'$. We can find $v'$ from~$v$ by tracing for each
intermediate ``tableau with empty square''~$\tilde T$ during the slide
$T\slid{T'}$ which of its entries~$i+1$ corresponds to the letter affected
by~$e_i$ in~$\wKan(\tilde T)$ (like before, we shall label this entry in
drawings by underlining it). This amounts to moving along with that entry if
it slides into another square, and switching to the entry to its right if an
entry~$i$ moves into its present column, giving rise to a
transition~$\smallsquares\Young(,0|\ul1,1)\to\Young(0,|1,\ul1)\,$.

We now distinguish two cases, depending on whether or not the path of the
inward slide applied to~$T$ (i.e., the set of squares whose entries move) is
the same as for the slide applied to~$e_i(T)$. Suppose first that the paths
are the same, then the entry of~$v$ makes the same move, if any, during the
two slides; since this entry is~$i+1$ in one case and $i$ in the other, it
must be unique in its column both before and after the slide. By the above
description of~$v'$, we see that it coincides with final position of the entry
of~$v$; this implies the theorem for this first case. Suppose next that the
paths of the two slides differ. This means that $v$ occurs in exactly one of
these paths; since the rule for an inward jeu de taquin slide is to select the
smallest candidate entry to move, $v$~must be in the path of the slide applied
to~$e_i(T)$, where it has entry~$i$ rather than~$i+1$. We can see as follows
that this move of the entry~$i$ of~$v$ in~$e_i(T)$ cannot be to the left: this
would mean that in the slide $T\slid T'$, where the entry~$i+1$ of~$v$ does
not move, the square to its left is filled by another entry~$i+1$ sliding up,
but then we find $v'=v$, which is absurd since $e_i(T')$ must have weakly
increasing rows. Therefore the move of the entry~$i$ of~$v$ in~$e_i(T)$ must
be upwards, and it is the first move of the slide. It is a transition
$\smallsquares\Young(,0|0,1)\to\Young(0,0|,1)$ whereas in~$T$ we have
$\smallsquares\Young(,0|\ul1,1)\to\Young(0,|1,\ul1)\,$ (the presence of the
entry at the bottom right follows from the fact that $v$ is the variable
square: the Kanji reading $\wKan(R)$ of the subtableau~$R$ of~$T$ consisting
of the columns to the right of~$v$ is anti-dominant for~$i$). The latter
transition may be followed by a number of similar transitions further to the
right, where each time the anti-dominance of~$\wKan(R)$ guarantees the
presence of the entry at the bottom right; eventually they must be followed by
a final transition $\smallsquares\Young(0,|1,\ul1)\to\Young(0,\ul1|1,)$. The
relevant part of the slide for~$T$ looks like
$\smallsquares\Young(,0,0,0|\ul1,1,1,1)\slid\Young(0,0,0,\ul1|1,1,1,)\,$, and
the location of the underline after the slide indicates the variable
square~$v'$ in~$T'$. The corresponding part of the slide for~$e_i(T)$ is
$\smallsquares\Young(,0,0,0|0,1,1,1)\slid\Young(0,0,0,0|1,1,1,)\,$, whose
result matches~$e_i(T')$; this proves the theorem for this second case.
\QED

Figure~$1$ illustrates the commutation of jeu de taquin and coplactic
operations. In this example the operation~$e_i$ with smallest possible~$i$ is
always chosen; several of them have been combined at each step, to prevent the
display from getting excessively large. The reader is urged to study the
transitions carefully.
\midinsert
$$
\ifnarrow \let\quad=\enspace \let\medskipamount=\smallskipamount \fi
\matrix
{
 \Young(0,0,0,0,0|1,1,1,1,1|2,2,2,2|3|4)
&\Skew(1:0,0,0,0|0:0,1,1,1,1|0:1,2,2,2|0:2,3|0:4)
&\Skew(1:0,0,0,0|1:1,1,1,1|0:0,2,2,2|0:1,3|0:2,4)
&\Skew(2:0,0,0|1:0,1,1,1|0:0,1,2,2|0:1,2,3|0:2,4)
&\Skew(3:0,0|1:0,0,1,1|0:0,1,1,2|0:1,2,2,3|0:2,4)
\cr\noalign{\medskip}
 \Young(0,0,0,0,0|1,1,1,1,1|2,2,2,5|3|5)
&\Skew(1:0,0,0,0|0:0,1,1,1,1|0:1,2,2,5|0:2,3|0:5)
&\Skew(1:0,0,0,0|1:1,1,1,1|0:0,2,2,5|0:1,3|0:2,5)
&\Skew(2:0,0,0|1:0,1,1,1|0:0,1,2,5|0:1,2,3|0:2,5)
&\Skew(3:0,0|1:0,0,1,1|0:0,1,1,2|0:1,2,3,5|0:2,5)
\cr\noalign{\medskip}
 \Young(0,0,0,0,0|1,1,1,1,1|2,4,4,5|3|5)
&\Skew(1:0,0,0,0|0:0,1,1,1,1|0:1,2,4,5|0:3,4|0:5)
&\Skew(1:0,0,0,0|1:1,1,1,1|0:0,2,4,5|0:1,4|0:3,5)
&\Skew(2:0,0,0|1:0,1,1,1|0:0,1,2,5|0:1,4,4|0:3,5)
&\Skew(3:0,0|1:0,0,1,1|0:0,1,1,2|0:1,4,4,5|0:3,5)
\cr\noalign{\medskip}
 \Young(0,0,0,0,0|1,1,1,3,3|2,4,4,5|3|5)
&\Skew(1:0,0,0,0|0:0,1,1,3,3|0:1,2,4,5|0:3,4|0:5)
&\Skew(1:0,0,0,0|1:1,1,3,3|0:0,2,4,5|0:1,4|0:3,5)
&\Skew(2:0,0,0|1:0,1,3,3|0:0,1,2,5|0:1,4,4|0:3,5)
&\Skew(3:0,0|1:0,0,1,3|0:0,1,2,3|0:1,4,4,5|0:3,5)
\cr\noalign{\medskip}
 \Young(0,0,0,1,1|1,1,2,3,3|2,4,4,5|3|5)
&\Skew(1:0,0,1,1|0:0,1,2,3,3|0:1,2,4,5|0:3,4|0:5)
&\Skew(1:0,0,1,1|1:1,2,3,3|0:0,2,4,5|0:1,4|0:3,5)
&\Skew(2:0,1,1|1:0,1,3,3|0:0,2,2,5|0:1,4,4|0:3,5)
&\Skew(3:0,1|1:0,1,1,3|0:0,2,2,3|0:1,4,4,5|0:3,5)
\cr
}
$$
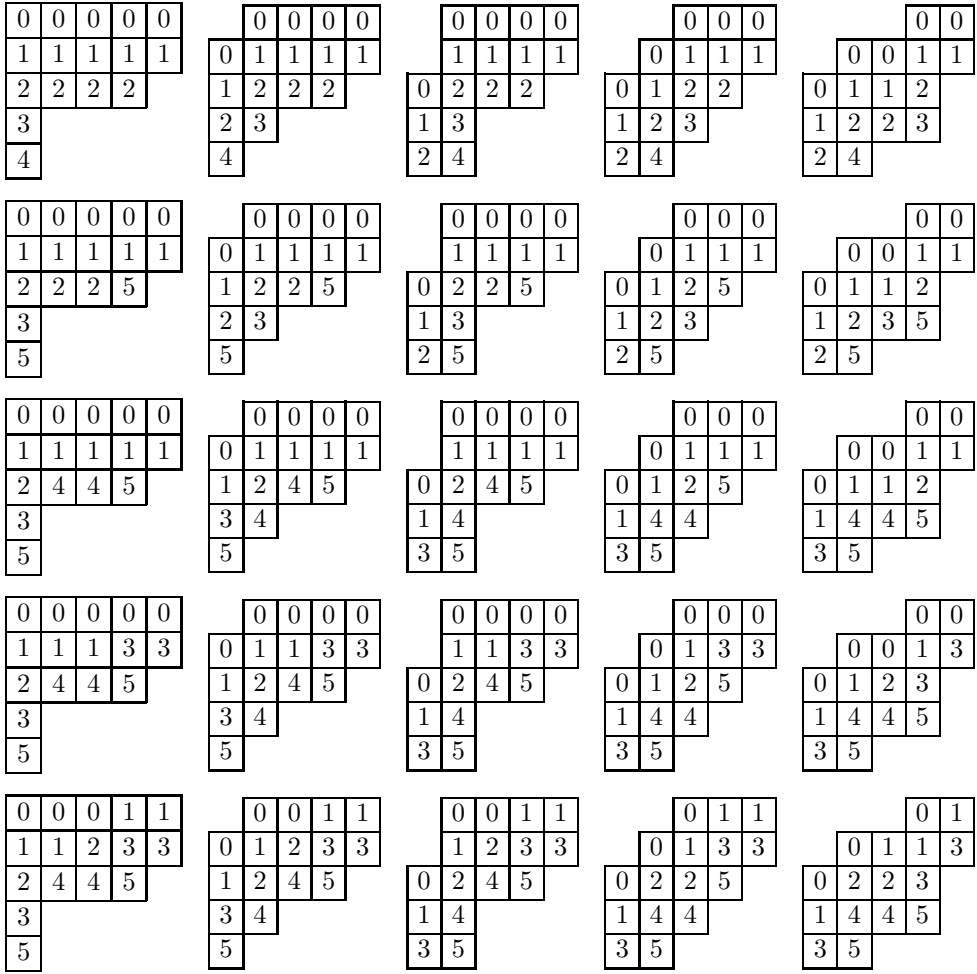
\centerline{Figure 1. Commutation of jeu de taquin (right to left) and raising
operations (upwards).}
\centerline{Raising operations are grouped: $e_0,e_0,e_1$; \
$e_2,e_1,e_2,e_1$; \ $e_3,e_2,e_3,e_2$; \ $e_4,e_4,e_3,e_2$.}
\endinsert

\proclaim Corollary.
Equation (\LRfibereq) and hence the Littlewood-Richardson rule have been
established.

\proof
Using corollary~\LRdomincorr, theorem~\copcommthm\ implies that jeu de taquin
preserves the property of being a Littlewood-Richardson tableau.
For~$L\in\Tab(\chi)$ there exists a tableau of partition shape
$P\in\Tab(\nu/0)$ with $L\slid P$; then $P$ is a Littlewood-Richardson tableau
if and only if~$P=\1_\nu$, in which case $\nu=\wt L$.
\QED

Note that we did not use the commutation statement of theorem~\copcommthm,
just the (easier) statement about when $e_i(T)$ is defined. Neither did we use
confluence of jeu de taquin here, although it was used to obtain
corollary~\altLRR. But in fact theorem~\copcommthm\ independently proves this
confluence, and more.

\proclaim Corollary.
On the set of all skew semistandard tableaux, the two rewrite systems defined
by inward jeu de taquin slides, respectively by the raising operations~$e_i$,
are both confluent (their normal forms are the semistandard Young tableaux,
and the Littlewood-Richardson tableaux, respectively). Moreover the two normal
forms of a skew semistandard tableau~$T$ uniquely determine~$T$.

\proof
Since the two rewrite systems commute in the precise sense of
theorem~\copcommthm, the set of normal forms for one system is closed under
the rewrite rules of the other system. Confluence of that rewrite system on
this set of normal forms is clear, since the only tableaux that are normal
forms for both systems simultaneously are the tableaux~$\1_\\$
for~$\\\in\Part$, and for a tableau that is a normal form for one of the
systems the pertinent value of~$\\$ can be directly read off as the shape (in
the case of semistandard Young tableaux) respectively as the weight (in the
case of Littlewood-Richardson tableaux).

Now let $T$ be any skew semistandard tableau, let $P$ be a tableau of
partition shape obtained from~$T$ by a sequence of inward jeu de taquin
slides, and let $\tilde P$ designate the sequence of shapes of the
intermediate tableaux. Similarly let $L$ be a Littlewood-Richardson tableau
obtained from~$T$ by applying a sequence~$\tilde L$ of raising operations. As
normal forms for one system, $P$~and~$L$ each have a unique normal form for
the other system, and by theorem~\copcommthm, one way to reach that normal
form is to apply the sequence~$\tilde L$ of raising operations to~$P$,
respectively to apply to~$L$ a sequence of jeu de taquin slides with $\tilde
P$ as sequence of intermediate shapes; moreover the two normal forms are the
same tableau~$\1_\\$ (where $P\in\Tab(\\/0)$ and~$\wt L=\\$). Now $P$ can be
reconstructed from~$\1_\\$ by reversing the the raising operations of~$\tilde
L$, and is therefore independent of the sequence~$\tilde P$ used to find it;
similarly $L$ can be reconstructed using~$\tilde P$, and therefore independent
of~$\tilde L$. If $L$~and~$P$ are given, then $T$ can be reconstructed by a
similar process: for instance, a sequence~$\tilde P$ of shapes can be found by
reducing~$L$ by inward jeu de taquin slides into~$\1_\\$, and since the latter
is also obtainable from~$P$ by raising operations, a sequence of outward jeu
de taquin slides can be applied to~$P$ involving the shapes of~$\tilde P$ in
reverse order, which produces~$T$ as result.
\QED

It may be noted that parts of the argument above could have been formulated in
more concrete terms, e.g., by using a skew standard tableau to encode the
information of~$\tilde P$ and using tableau switching; we have not done this
in order to stress the conceptual simplicity and symmetry with respect to the
two rewrite systems, and the fact that no detailed knowledge about these
systems is used. On the other hand the computation of~$T$ from $L$~and~$P$
above can easily be seen to be equivalent to that of~$\phi(L,P)$ in
corollary~\LRjdtcor, so that the following defines a specialisation
of~$\Rob'$, as claimed at the end of~\Sec2.

\proclaim Corollary. \Robbijdescr
For $\chi\in\Sk$ and $n\in\N$, Robinson's bijection~$\Rob$ (i.e., one
satisfying the specification~(\Robinson)) can be defined by
$\Rob(T)=(\Rob_0(T),\Rob_1(T))$ for $T\in\Tab(\chi;n)$, where
$\Rob_0(T)\in\LR(\chi,\nu)$ is the normal form of~$T$ for the rewrite system
of raising operations, and the $\Rob_1(T)\in\Tab(\nu;n)$ is the normal form
of~$T$ for the rewrite system of inward jeu de taquin slides.
\QED

\subsection Jeu de taquin on companion tableaux.

Having completed our discussion of the Littlewood-Richardson rule proper, it
is worth while to point out a remarkable connection between coplactic
operations and jeu de taquin performed on companion tableaux. Note that we
have already found in corollary~\LRdomincorr\ that $T$ has a companion tableau
that is a normal form for jeu de taquin (i.e., $T$~is a Littlewood-Richardson
tableau) if and only if $T$ itself is a normal form for raising operations.
The connection we shall establish generalises this; it is closely related to
the theory of pictures, where instead of coplactic operations one has a second
form of jeu de taquin, with a similar commutation theorem \ref{pictures,
theorem~5.3.1}. The symmetry exhibited by this connection has some important
implications, which we shall only indicate briefly here. Firstly, it implies
the fundamental symmetry of the Robinson-Schensted correspondence (see for
instance \ref{eljc, 3.2}), which is more manifest for Schensted's formulation
of the correspondence \ref{Schensted} than for Robinson's bijection; however,
pictures clearly exhibit the link between these formulations \ref{pictures,
theorem~5.2.3}. Secondly, the connection clarifies the bijection
corresponding to~$c_{\\,\mu}^\nu=c_{\mu,\\}^\nu$: as remarked
at the end of~\Sec1, this bijection involves tableau switching performed on a
companion tableau of (part of) the Littlewood-Richardson tableau. A detailed
discussion of this bijection, and a generalisation, can be found in \ref{path
jeu de taquin}.

\proclaim Theorem. \jdtcoplink
Let $\bT\in\Tab(\nu/\kappa)$ be a companion tableau of~$T\in\Tab(\chi)$, and
let $\bT'\in\Tab(\nu'/\kappa')$ be obtained from~$\bT$ by an inward jeu de
taquin slide into a square in row~$k$, ending in row~$l$. Then the sequence
$T=T_k,T_{k+1},\ldots,T_l\in\Tab(\chi)$ given by $T_{i+1}=e_i(T_i)$ for $k\leq
i<l$ is well defined, and $\bT'$ is a companion tableau of~$T_l$; moreover
$T_l$ is the first tableau in the sequence that is $\kappa'$-dominant.

\proof
We augment the word~$\wSem(T)$ obtained from~$T$ by the Semitic reading order
with additional data, by attaching a subscript called ``ordinate'' to each
letter: for each~$i\in\N$, the ordinates of the letters~$i$ of~$\wSem(T)$
increase from left to right by unit steps, starting at~$\kappa_i$. We thus
obtain an augmented word~$x$ that describes the standardisation of~$\bT$: the
sequence of squares~$(i,j)$ added to~$Y(\kappa)$ corresponds to the sequence
of letters~$i_j$ (letter~$i$ with ordinate~$j$) of~$x$ from left to right.
During the application of an inward slide to the standardisation of~$\bT$,
transforming it into that of~$\bT'$, we modify~$x$ correspondingly: for each
move of an entry from a square~$(i,j+1)$ or~$(i+1,j)$ to $(i,j)$ we replace
the letter~$i_{j+1}$ respectively $(i+1)_j$ in~$x$ by~$i_j$. Due to these
changes, the underlying word of~$x$ (obtained by forgetting the ordinates)
will be transformed from~$\wSem(T_k)$ into~$\wSem(T_l)$: replacement
of~$i_{j+1}$ by~$i_j$ does not change the underlying word, while we claim that
changing $u(i+1)_jv$ into~$ui_jv$ corresponds to an application of~$e_i$ to
the underlying word, in other words that the underlying word of $u$ is
anti-dominant for~$i$ and that of~$v$ dominant. From the correspondence with
an intermediate tableau for the slide, it follows that for any letter~$i_{j'}$
occurring in~$u$ one has $j'<j$, and there is a letter $(i+1)_{j'}$ to the
right of it in~$u$, and similarly that $j'>j$ for any letter~$(i+1)_{j'}$
in~$v$, and there is a letter~$i_{j'}$ to its left; from this our claim is
easily deduced. By construction the standardisation of~$\bT'$ is described by
the word into which $x$ is finally transformed, whose underlying word
is~$\wSem(T_l)$; since $\wt\bT'=\wt\bT$, this implies that $\bT'$ is a
companion tableau of~$T_l$. It remains to show that $T_i$ is not
$\kappa'$-dominant for $i<l$: apply the test for $\kappa'$-dominance of
proposition~\LRlatticechar\ to~$T_i$, and consider the value of~$\alpha$ just
before the entry is visited that corresponds to the letter~$(i+1)_j$ in~$x$
that is changed into~$i_j$; the fact that that change takes place means that
$\alpha_i=\alpha_{i+1}=j$, whence the test will fail.
\QED

As an example, we apply this construction to the pair~$T,\bT$ of companion
tableaux of~(\companionex), and~$k=0$. The word
$\let~=\,\wSem(T)=1~0~3~1~1~0~3~2~2~0~5~4~4~1~5~3$, gets augmented as
$\let~=\,1_4~0_6~3_2~1_5~1_6~0_7~3_3~2_4~2_5~0_8~5_0~4_1~4_2~1_7~5_1~3_4$,
since $\kappa=(6,4,4,2,1,0)$.
The slide into the square~$(0,5)$ causes the changes $0_5\gets0_6$, \
$0_6\gets1_6$, and $1_6\gets1_7$, resulting in the augmented word
$\let~=\,1_4~0_5~3_2~1_5~0_6~0_7~3_3~2_4~2_5~0_8~5_0~4_1~4_2~1_6~5_1~3_4$; on
the level of the underlying word this amounts to application of~$e_0$, and the
pair of companion tableaux
$$
  e_0(T)=\Skew(3:0,1|1:0,0,1,3|0:0,2,2,3|0:1,4,4,5|0:3,5),
\quad
  \bT'=\Skew(5:0,1,1,2|4:0,1,3|4:2,2|2:1,2,4|1:3,3|0:3,4),
\nn
$$
shows that $e_0(T)$ is $\kappa'$-dominant for~$\kappa'=(5,4,4,2,1,0)$. A
further slide into the square~$(0,4)$ causes the changes $0_4\gets1_4$, \
$1_4\gets1_5$, and $1_5\gets2_5$, resulting in
$\let~=\,0_4~0_5~3_2~1_4~0_6~0_7~3_3~2_4~1_5~0_8~5_0~4_1~4_2~1_6~5_1~3_4$, and
$$
  (e_1\after e_0)(e_0(T))=\Skew(3:0,0|1:0,0,1,3|0:0,1,2,3|0:1,4,4,5|0:3,5),
\quad
  \bT''=\Skew(4:0,0,1,1,2|4:1,2,3|4:2|2:1,2,4|1:3,3|0:3,4).
\nn
$$
Writing the letters of the augmented words back into~$T$, one finds the
transitions
$$
  \Skew(3:            0_6,1_4
       |1:    0_7,1_6,1_5,3_2
       |0:0_8,2_5,2_4,3_3
       |0:1_7,4_2,4_1,5_0
       |0:3_4,5_1)
\quad\buildrel \textstyle e_0 \over \longrightarrow\quad
  \Skew(3:            0_5,1_4
       |1:    0_7,0_6,1_5,3_2
       |0:0_8,2_5,2_4,3_3
       |0:1_6,4_2,4_1,5_0
       |0:3_4,5_1)
\quad\buildrel \textstyle e_1\after e_0 \over \longrightarrow\quad
  \Skew(3:            0_5,0_4
       |1:    0_7,0_6,1_4,3_2
       |0:0_8,1_5,2_4,3_3
       |0:1_6,4_2,4_1,5_0
       |0:3_4,5_1).
\label(\imglisex)
$$
This representation provides another perspective on proposition~\coptabprop.
First observe that the ordinates of the entries of~$T$ with a fixed value~$i$
increase \emph{from right to left} by unit steps, starting with the
ordinate~$\kappa_i$. It follows that reading the augmented entries according
to any valid reading order~`$\leqr$' will produce a properly augmented word;
stated differently, augmenting (for~$\kappa$) $w_r(T)$ and then writing the
augmented letters back into~$T$ always produces the same augmented tableau.
Also, the values of the entries of the augmented tableau~$T$ with a fixed
ordinate~$j$ increase from top to bottom, by what was said in~\Sec\readingsec.
The fact that during the first transition we made the
replacement~$0_5\gets0_6$ rather than $0_5\gets1_5$ is due to the fact that
$0_6$ precedes~$1_5$ in the augmented word~$\wSem(T)$, and similarly the
replacement~$0_6\gets1_6$ was made because $1_6$ precedes~$0_7$ (for the final
replacement~$1_6\gets1_7$ there is no alternative). As can be seen, these
ordering relations are unchanged when an augmented word~$w_r(T)$ is used
instead of~$\wSem(T)$, due to the placement of these entries in the augmented
tableau~$T$, and the same is true for the relevant pairs of entries
$(1_4,0_5)$, $(1_5,2_4)$ and~$(2_5,1_6)$ in the second augmented tableau
in~(\imglisex). The fairly easily proved fact that this is always the case
(cf.\ \ref{pictures, lemma~5.1.2}) provides alternative proof of
proposition~\coptabprop.

Using theorem~\jdtcoplink, one can express coplactic operations on~$T$ in
terms of jeu de taquin slides on a suitable companion tableau of~$T$: to
compute $e_k(T)$ (when defined), it suffices to find a companion tableau~$\bT$
for which one has $l=k+1$ in the theorem, and this can be achieved by choosing
$\nu/\kappa$ in such a way that $\kappa_k-\kappa_{k+1}$ is minimal, and that
$\nu_{k+2}<\kappa_k$ (the latter condition can be weakened considerably).
Conversely, if one wishes to express jeu de taquin slides on a companion
tableau of~$T$ in terms of coplactic operations on~$T$, it is necessary to
replace the condition of $\kappa'$-dominance in the theorem by a condition
stated in terms of coplactic graphs; this is possible using the following
proposition.

\proclaim Proposition. \jdttocop
Let~$T\in\Tab(\chi,n)$, and let $\nu/\kappa$ be a skew shape with
$\nu-\kappa=\wt T$. Then the following are equivalent
\eqitem $T$ is $\nu/\kappa$-dominant;
\eqitem for all~$i\in\upto{n-1}$ the operation~$e_i$ cannot be applied more
        than $\kappa_i-\kappa_{i+1}$ times successively to~$T$;
\eqitem for all~$i\in\upto{n-1}$ the operation~$f_i$ cannot be applied more
        than $\nu_i-\nu_{i+1}$ times successively to~$T$.

For instance, for the tableau~$T$ of~(\companionex) used in the example above
we have $\kappa=(6,4,4,2,1,0)$ and $\nu=(9,8,6,5,3,2)$, so the proposition
states that $e_1$ cannot be applied to~$T$, that $e_3$, $e_4$, $f_0$, $f_2$
and~$f_4$ cannot be applied more than once to~$T$, and that $e_0$, $e_2$,
$f_1$, and~$f_3$ cannot be applied more than twice; this can be verified, and
one finds moreover that the statement can be sharpened for~$e_3$ (which cannot
be applied) and for~$f_3$ (which can be applied only once), which corresponds
to the fact that~$T$ is also $\nu'/\kappa'$-dominant, with
$\nu'=(8,7,5,4,3,2)$ and $\kappa'=(5,3,3,1,1,0)$.

\proof
Let $T*\Can\kappa$ denote a tableau of a shape representing~$\chi*(\kappa/0)$
whose restriction to the factor~$\chi$ of~$\chi*(\kappa/0)$ is a translation
of~$T$, while its restriction to the complementary factor~$\kappa/0$ is a
translation of~$\Can\kappa$. Then \eqitemnr1 is equivalent to $T*\Can\kappa$
being $0$-dominant, and hence to
$\wSem(T*\Can\kappa)=\wSem(\Can\kappa)\wSem(T)$ being dominant; also the
number of times $e_i$ or~$f_i$ can be applied to~$T$ is the same as
for~$\wSem(T)$. Fixing any $i\in\upto{n-1}$, set $c=\kappa_i-\kappa_{i+1}$ and
$d=\nu_i-\nu_{i+1}$; we have to show that for a word~$w\in\upto{n}^*$ of
weight~$\nu-\kappa$ to which $e_i$ and~$f_i$ can be successively applied
exactly $r$ respectively $s$~times, there is equivalence between: \
(1)~$\wSem(\Can\kappa)w$ is dominant for~$i$, \ (2)~$r\leq c$, \ (3)~$s\leq
d$. Since $\wt w=\nu-\kappa$ implies $r-s=c-d$, (2)~and~(3) are equivalent; we
shall prove equivalence of (1)~and~(2), dropping the assumption
$\wt{w}=\nu-\kappa$, since $\nu$ is no longer relevant. Subwords that are
neutral for~$i$ neither affect condition~(1) nor the value of~$r$; removal of
such subwords reduces~$w$ to~$(i+1)^ri^s$ (where exponentiation signifies
repetition of letters), and $\wSem(\Can\kappa)$ to~$i^c$. It remains to show
that $i^c(i+1)^ri^s$ is dominant for~$i$ if and only if~$r\leq c$, which is
obvious.
\QED

\proclaim Corollary.
Let $T\in\Tab(\chi)$ and $T'\in\Tab(\chi')$ be jeu de taquin equivalent, and
$\nu/\kappa$ a skew shape. Then $T$ is $\nu/\kappa$-dominant if and only if\/
$T'$ is so, in which case the companion tableaux of\/ $T$~and~$T'$ of shape
$\nu/\kappa$ are dual equivalent.

\proof
By theorem~\copcommthm, $T$ and~$T'$ have isomorphic coplactic graphs, and
proposition~\jdttocop\ then shows that $T$ is $\nu/\kappa$-dominant if and
only if $T'$ is. Let $\bT$ and~$\bT'$ respectively be the companion tableaux
of $T$~and~$T'$ of shape $\nu/\kappa$, and consider an inward jeu de taquin
slide into the same square applied to each of them. By theorem~\jdtcoplink,
the result will in either case be a companion tableau of a tableau obtained by
a sequence of raising operations from $T$ respectively from~$T'$. Moreover,
since the condition determining the length of those sequences can be expressed
in terms of the coplactic graph by proposition~\jdttocop, the two sequences
will be identical. Therefore the jeu de taquin slides applied to $\bT$
and~$\bT'$ result in tableaux of equal shape that are companion tableaux of
jeu de taquin equivalent tableaux; the same will hold after further inward
slides. Proposition~\revdeq\ shows that this suffices to establish dual
equivalence.
\QED

We have shown that jeu de taquin slides performed on~$T$ commute with jeu de
taquin slides performed on a companion tableau of~$T$; this is essentially
\ref{pictures, theorem~5.3.1}. One may ask how practical this alternative
method of computing coplactic operations is. For single coplactic operations,
or following prescribed paths in the coplactic graph (such as the one used in
figure~1), the method is cumbersome, as one (repeatedly) has to adapt the
shape~$\nu/\kappa$ to single out the right coplactic operation. (Note however
that the path that always applies $e_i$ with the \emph{largest} possible~$i$
can be found without such adaptations.) On the other hand, when the goal is
just to compute~$\Rob_0(T)$, the fact that one often gets more than one
coplactic operation at a time is an advantage, as is the absence of a need to
repeatedly recompute globally defined quantities, such as the location of the
square affected by some~$e_i$. In addition, this method is less error prone
when used for manual computation. For instance in~(\companionex), by computing
$$
  {\maybesmall\Skew(6:0,1,2|4:0,1,1,3|4:2,2|2:1,2,4|1:3,3|0:3,4)}
\slid
  {\maybesmall\Young(0,0,1,1,2|1,1,2,2,3|2,3,3,4|3|4)},
\hbox{ one finds that }
 \Rob_0\left(\,{\maybesmall\Skew(3:0,1|1:0,1,1,3|0:0,2,2,3|0:1,4,4,5|0:3,5)}
     \,\right)
=
 {\maybesmall\Skew(3:0,0|1:0,0,1,1|0:0,1,1,2|0:1,2,2,3|0:2,4)}.
$$
\endofchapter
 \par

\ifx\macrosloaded\undefined
  \input marxmax
  

\newif\ifnarrow
\ifdim \hsize<14.5cm \narrowtrue \else\narrowfalse \fi

\def\maybesmall{\ifnarrow\smallsquares\fi}

\def\ol#1
{\vbox{\vskip0.3pt \hrule \vskip 0.8pt \hbox{$#1$}}}

\def\ul#1
{\vtop{\hbox{$#1$}\vskip 0.8pt \hrule \vskip0.3pt }}

\font\frak=eufm10
\font\csc=cmcsc10
\font\URLfont=cmtex10 
\hyphenchar\URLfont=-1 

\mathorddef\Part=P2 
\mathorddef\A=A0 \mathorddef\B=B0 
\def\GL{\boldmath{GL}}
\def\Ylat{(\Part,{\subset})}
\def\upto#1{[#1]}
\opdef\Pic
\altopdef\Tab{SST}
\opdef\ST
\opdef\LR
\opdef\JdT
\opdef\wt
\mathorddef\T=T2 
\mathorddef\Sk=S2 
\def\SFn{\Lambda_n}
\def\Can#1{\boldmath1_{#1}}
\def\Ylm{Y(\\/\mu)}
\def\bT{\bar T}
\def\LRt#1,#2{\LR(#1,#2)}
\def\JdTt#1,#2{\Tab(#1)^{\slid{#2}}}
\mathorddef\Rob=R2 
\def\StT#1{\ST(#1)}
\def\snu{^{\slid\nu}}
\def\Uqgln{U_q(\hbox{\frak gl}_n)}
\def\leqr{\leq_r}
\def\leqKan{\leq_{\romath K}}

\def\wKan{w_{\romath K}}
\def\wSem{w_{\romath S}}

\mathdef\slid=Rel \triangleright
\mathdef\dils=Rel \triangleleft

\long\def\quotation#1
{{\smallskip\narrower\it\noindent\llap{``\/}#1\unskip''\smallskip}\noindent}

{\catcode`\.=\active\catcode`\/=\active\catcode`\_=\active
 \gdef\URL
 {\begingroup\def~{\char`\~}%
  \catcode`\.=\active\catcode`\/=\active\catcode`\_=\active
  \let_=\_%
  \def.{\char`\.\penalty\hyphenpenalty}\def/{\char`\/\penalty\hyphenpenalty}%
  \URLfont\URLnext
}}
\def\URLnext#1{http://#1\endgroup}

\hyphenation
{semi-standard Yamano-uchi mono-mial mono-mials inter-related combi-na-torics
 combi-na-torial algo-rithm valid-ity}

  \writelab{}
  \def\jobname{lrr}
\readrefs
  \let\endofchapter=\bye
  \input lrr1.lab
  \input lrr2.lab
  \input lrr3.lab
\fi

\newsection Some historical comments.

Now that we have seen the Littlewood-Richardson rule from a modern
perspective, let us look at some elements of its intriguing history, in
particular the two papers \ref{Littlewood Richardson}~and~\ref{Robinson},
written in the 1930's. 

\subsection The paper by Littlewood and Richardson.

The paper in which the Littlewood-Richardson rule is first stated, is mainly
concerned with symmetric group characters and Schur functions (a term
introduced in that paper, though mostly contracted to ``S-functions''). Out of
16~sections, only \Sec8 deals with the multiplication of S-functions.
Remarkably, semistandard tableaux do not occur explicitly, and in particular
are not used in the definition of Schur functions; instead these are expressed
in terms of power sum symmetric functions using symmetric group characters. In
fact no attempt is made at all to express Schur functions in terms of
monomials, or even of symmetric monomials~$m_\\$, although curiously the
opposite \emph{is} done,~in~\Sec5. Semistandard tableaux do occur implicitly,
as follows. For~$\mu=(r)\in\Part_r$ one has $s_\mu=h_r$, and every skew
semistandard tableau of weight~$\mu$ is automatically a Littlewood-Richardson
tableau; call the shape of such a tableau a horizontal strip. As starting
point for multiplication of S-functions, it is shown that $s_\\h_r$ is the sum
of all $s_\nu$ with $\nu/\\$ a horizontal strip (this agrees with our
proposition~\LRswitchprop). By iteration this implies that
$s_\\\fold{}(h\alpha_0..\alpha_l)$ is the sum of $s_\nu$ taken over all ways
to successively add horizontal strips of sizes $\li(\alpha0..l)$ to~$Y(\\)$
yielding~$Y(\nu)$. The number of ways to so obtain a given partition~$\nu$ is
$\Card{\setof{T}\in\Tab(\nu/\\):\wt{T}=\alpha\endset}$. For the general case
of multiplying S-functions, the following rule is formulated (\ref{Littlewood
Richardson,~p.\thinspace119}):
\quotation
{{\csc Theorem III.}---Corresponding to two $\romath S$-functions
 $\{\li(\\1..p)\}$,~$\{\li(\mu1..q)\}$ build tableaux $\A$~and~$\B$ as in
 Theorem~II. Then in the product of these two functions, the coefficient of
 any $\romath S$-function $\{\nu_1,\nu_2,\ldots\}$ is equal to the number of
 compound tableaux including all of the symbols of\/ $\A$~and~$\B$, and
 corresponding to~$\{\nu_1,\nu_2,\ldots\}$, that can be built according to the
 following rules.

 Take the tableau~$\A$ intact, and add to it the symbols of the first row
 of\/~$\B$. These may be added to one row of\/~$\A$, or the symbols may be
 divided without disturbing their order, into any number of sets, the first
 set being added to one row of\/~$\A$, the second set to a subsequent row, the
 third to a row subsequent to this, and so on. After the addition no row must
 contain more symbols than a preceding row, and no two of the added symbols
 may be in the same column.

 Next add the second row of symbols from~$\B$, according to the same rules,
 with this added restriction. Each symbol from the second row of\/~$\B$ must
 appear in a later row of the compound tableau than the symbol from the first
 row in the same column.

 Similarly add each subsequent row of symbols from~$\B$, each symbol being
 placed in a later row of the compound tableau than the symbol in the same
 column from the preceding row of\/~$\B$, until all the symbols of\/~$\B$ have
 been used.
}
The tableaux~$\A$ and~$\B$ are Young diagrams of shapes $\\$ and~$\mu$, filled
with formal symbols. The rows of~$\B$ are rearranged into horizontal strips
filling~$Y(\nu/\lambda)$; as indicated above such a rearrangement corresponds
to some $T\in\Tab(\nu/\\)$ with~$\wt T=\mu$. The restriction added in the last
two paragraphs depends on the particular rule given to ensure that each
horizontal strip is obtained in only one way from a given row of~$\B$, namely
that within each row of~$\B$ the number of the destination row of each symbol
increases weakly from left to right. If in~$\B$ one labels each symbol with
this number, then the restriction says that that these numbers increase
strictly down columns. The numbers then form a tableau~$\bT\in\Tab(\mu/0)$,
that can be seen to be a companion tableau of~$T$ in the sense of the current
paper. In fact if we interpret the rearrangement as a bijection~$Y(\mu)\to
Y(\nu/\lambda)$, then we almost arrive at the notion of pictures; the only
difference is that when some symbols from a row of~$\B$ are moved to a common
row, this happens ``without disturbing their order'', rather than by reversing
their order, as in pictures.

Although claimed to be generally valid, Theorem~III is only proved
in~\ref{Littlewood Richardson} when $\mu$ has at most two parts (the first
sentence after the statement of the Theorem is ``No simple proof has been
found that will demonstrate it in the general case''; indeed ``simple'' can
be omitted). We already mentioned the particularly simple case $\mu=(r)$. The
proof for a partition $\mu=(q,r)$ with two parts $q\geq r>0$ is based on the
identity $s_\mu=h_qh_r-h_{q+1}h_{r-1}$ (an instance of a determinantal formula
known as the Jacobi-Trudi identity), which means that $c^\nu_{\lambda,\mu}$
can be found by subtracting the coefficient of~$s_\nu$ in~$s_\\h_{q+1}h_{r-1}$
from the one in~$s_\\h_qh_r$; these coefficients can be determined by counting
the tableaux in~$\Tab(\nu/\\,2)$ of weights $(q,r)$ and~$(q+1,r-1)$,
respectively. To show that the difference matches the number of tableaux given
by Theorem~III, the rule is reformulated in terms of lattice permutations
(cf.~proposition~\LRswitchprop), and a bijective correspondence
in~$\Tab(\nu/\\,2)$ is given (in fact the one defined by our~$e_0$) between
tableaux of weight~$(q,r)$ that do not correspond to lattice permutations, and
tableaux of weight~$(q+1,r-1)$; it is verified that the transformation
preserves semistandardness.

It is remarkable that the authors state their rule as a Theorem when, by their
own admission, they only have a proof for some very simple cases. In part this
can be attributed to the general attitude at the time, which appears to have
been that combinatorial statements are less in need of a proof than, say,
algebraic statements; on the other hand, they do devote three pages to a proof
of the special cases. It appears that they viewed their rule mainly as a
computational device, useful to find a result, whose correctness may then be
verified by other means. They work out a complete example for the computation
of $s_{(4,3,1)}s_{(2,2,1)}$ (in their notation $\{4,3,1\}\times\{2^21\}$),
displaying all~$34$ tableaux contributing to the result. They use
$$
\let\quad=\enspace
  \pmatrix{a,&b,&c,&d\cr e,&f,&g\hfill\cr h\hfill\cr},
\quad
  \pmatrix{\alpha,&\beta\cr\gamma,&\delta\cr\eps\hfill\cr},
$$
as tableaux~$\A$, $\B$; replacing the symbols of~$\A$ by $0$'s as they
don't move, they then display tableaux like
$$\let\quad=\enspace
  \matrix{0&0&0&0&\alpha\cr0&0&0&\beta\cr0&\gamma&\delta\cr\eps\cr};
$$
note the striking resemblance with pictures as displayed in~(\pictureex) (but
the order of $\gamma$~and~$\delta$ is unchanged). It is clear that they use a
geometric criterion ($\gamma$ should remain below~$\alpha$, \ $\delta$
below~$\beta$, and $\eps$ below~$\gamma$) rather than the lattice permutation
condition; their method is therefore an efficient one, as we discussed
following proposition~\LRswitchprop. The authors go on to indicate explicitly
how the resulting decomposition of the product can be verified by computing
the dimension of the corresponding (reducible) representation of the symmetric
group~$\Sym_{13}$. Interestingly, while their set of tableaux is correct, they
forget the contributions of three of them to the decomposition, whence the
dimension of the printed result, $398541$, falls short of the correctly
predicted dimension~$450450$, despite the claim that ``this equation proves
to be correct''.

One more curious point concerns the Theorem~II referred to in Theorem~III. It
describes the set of S-functions appearing in a product of two S-functions,
without giving their multiplicities. This would seem to contradict our
statement that no method is known to decide membership of that set, which does
not amount to finding (or failing to find) an appropriate
Littlewood-Richardson tableau. This is not so for two reasons. First, the
criterion given is even more impractical than using the Littlewood-Richardson
rule: it states that occurrence of $s_\nu$ in $s_\\s_\mu$ is equivalent to the
existence of a bijection $Y(\\*\mu)\to Y(\nu)$ mapping squares of any row to
distinct columns, and squares of any column to distinct rows. Enumerating this
set of bijections (or proving it empty) is certainly not easier than
enumerating its subset of pictures, which is equivalent to enumerating
$\LR(\nu/\mu,\\)$. Secondly the criterion is false; for instance the
bijection
$$
  \Skew(2:a,b|2:c,d|0:e,f|0:g,h)\quad\longrightarrow\quad
  \Young(a,b,e,h|c,d|f|g)
$$
satisfies the requirements, while $c_{(2,2),(2,2)}^{(4,2,1,1)}=0$. The proof
that is given in fact only estabishes the necessity of the condition, not its
sufficiency. Fortunately, Littlewood and Richardson will be remembered more
for a true Theorem they did not claim to prove, than for a false
Theorem they did claim to prove.

\subsection The paper by Robinson.

In~\ref{Robinson}, Robinson builds forth on these ideas, claiming to complete
the proof of the rule. The paper is quite difficult to read however, its
formulations extremely obscure, and its argumentation mostly implicit; we
shall now first try to summarise the argument as it appears to have been
intended. In a deviation from the previous paper, the proof does not use the
Jacobi-Trudi identity, but rather descending induction on the weight~$\mu$
with respect to~`$\prec$', based on the expression
$s_\mu=h_\mu-\sum_{\mu'\succ\mu}K_{\mu'\mu}s_{\mu'}$, where
$h_\mu=h_{\mu_0}h_{\mu_1}\cdots$, and the $K_{\mu'\mu}$ are non-negative
integer coefficients (nowadays called Kostka numbers). For instance, for
$\mu=(q,r)$, one has $s_\mu={h_qh_r-s_{(q+1,r-1)}-\cdots-s_{(q+r,0)}}$. The
correspondence~$\Rob$ of~(\Robinson) is defined, and used as follows to show
that each $L\in\LRt\nu/\\,\mu$ corresponds to an occurrence of~$s_\nu$
in~$s_\\s_\mu$. We already know that each $T\in\Tab(\nu/\\)$ with $\wt T=\mu$
corresponds to an occurrence of~$s_\nu$ in~$s_\\h_\mu$. If $\Rob_0(T)\neq T$,
then~$\Rob_0(T)$ corresponds by inductive assumption to an occurrence~$X$
of~$s_\nu$ in~$s_\\ s_{\mu'}$, where $\mu'=\wt\Rob_0(T)\succ\mu$. Viewing
$s_{\mu'}$ as a constituent of~$h_\mu$ and thereby $s_\\s_{\mu'}$ as a part
of~$s_\\h_\mu$, we let $T$ correspond to the occurrence of~$s_\nu$
in~$s_\\h_\mu$ matching~$X$. When $\Rob_0(T)=\Rob_0(T')$ for $T\neq T'$, then
the distinct tableaux $\Rob_1(T),\Rob_1(T')\in\Tab(\mu'/0)$ can be used to
label distinct constituents~$s_{\mu'}$ of~$h_\mu$; indeed there are
$K_{\mu'\mu}=\Card{\setof P\in\Tab(\mu'/0):\wt P=\mu\endset}$ such
constituents. The tableaux that remain, namely those with
$T=\Rob_0(T)\in\LRt\nu/\\,\mu$, must then correspond to the occurrences
of~$s_\nu$ in~$s_\\ s_\mu$, as claimed.

Now looking at what \ref{Robinson} actually says, we must first note that
neither the numbers $K_{\mu'\mu}$ nor semistandard Young tableaux occur
explicitly; rather the decomposition
$h_\mu=\sum_{\mu'\succeq\mu}K_{\mu'\mu}s_{\mu'}$ is given in the cryptic form
$h_\alpha=\sum\bigl[\prod S_{rs}^{\\_{rs}}\bigr](\alpha)$, with the following
``explanation'' (quoted literally from~\ref{Young}):
\quotation
{$S_{rs}$ {\rm where} $r<s$ represents the operation of moving one letter from
 the \hbox{$s$-th} row up to the \hbox{$r$-th} row, and the resulting term is
 regarded as zero, whenever any row becomes less than a row below~it, or when
 letters from the same row overlap,---as, for instance, happens when
 $\alpha_1=\alpha_2$ in the case of $S_{13}S_{23}$.
}
This is what is nowadays called Young's rule, although it is usually given in
a somewhat different from. The meaning of this original formulation appears to
be as follows. The summation is over certain collections $(\\_{rs})_{1\leq
r<s\leq n}$ with $\\_{rs}\in\N$, each of which gives rise to a monomial
$M=\prod S_{rs}^{\\_{rs}}$ in commuting indeterminates~$S_{rs}$. An operation
on Young diagrams is associated to such~$M$, where each factor~$S_{rs}$ moves
a square from row~$s$ up to row~$r$, all factors acting simultaneously; only
those monomials are considered in the summation for which this operation can
be applied to~$Y(\alpha)$, i.e., for each~$s$ one requires
$\sum_{r<s}\\_{rs}\leq\alpha_s$. According to our formula for~$h_\mu$, any
$P\in\Tab(\beta)$ with $\wt P=\alpha$ should correspond to such a monomial,
whose operation applied to~$Y(\alpha)$ results in $Y(\beta)$. This can be
achieved by taking $\\_{rs}=P_r^s$ (cf.~definition~\companiondef): then each
entry of~$P$ records the original row of its square (i.e., $P$ is obtained by
applying~$M$ as square-moving operation to~$\1_\alpha$, the entries following
their squares). It remains unclear how the phrase above can be interpreted as
describing only the monomials corresponding to such tableaux~$P$. The given
restrictions may be read as requiring $\beta$ to be a partition and forbidding
equal entries in the same column, but it is a mystery how the monotonicity of
columns is enforced: somehow $S_{12}S_{24}$ must be forbidden for
$\alpha=(1,1,1,1)$, while $S_{12}S_{23}S_{34}$ is allowed (they correspond to
$\smallsquares\smallsquares\Young(1,2|4|3)$
and~$\smallsquares\smallsquares\Young(1,2|3|4)$, respectively).

Robinson defines a process of transforming non-lattice permutations (i.e.,
words over an alphabet $\{c_1,c_2,\ldots\}$) into lattice permutations, which
he calls ``association~I''; incidentally, he attributes it to D.~E.
Littlewood. In our terminology it amounts to repeatedly applying the raising
operation~$e_i$ with smallest possible index~$i$ (where the letter~$c_j$ is
treated like the letter~$j$ in our description), until no further application
is possible; clearly, applying it to the Semitic reading~$\wSem(T)$ of a skew
tableau~$T$ of weight~$\alpha$, and writing the resulting lattice permutation
back into the shape of~$T$, this association corresponds to~$\Rob_0$ as
described in corollary~\Robbijdescr. Inspired by the fact that the monomials
$\prod S_{rs}^{\\_{rs}}$ are described as ``products of operations'', Robinson
associates such a monomial~$M$ to the sequence of operations~$e_i$ used in
determining association~I, by grouping together maximal sequences of
successive operations of the form $e_{s-1},\li(es-2..r+1),e_r$, replacing them
by~$S_{rs}$, and multiplying all of these. He calls the resulting
correspondence between the original word and this monomial ``association~II''.
For instance, the sequence used in figure~1 would be grouped as $(e_0)$,
$(e_0)$, $(e_1)$, $(e_2,e_1)$, $(e_2,e_1)$, $(e_3,e_2)$, $(e_3,e_2)$, $(e_4)$,
$(e_4,e_3,e_2)$, and would give rise to the monomial
$M=S_{01}^2S_{12}^{}S_{13}^2S_{24}^2S_{45}^{}S_{25}^{}$.

Robinson now states that~$M$ is one of the monomials that arises when Young's
rule is applied for~$h_\alpha$ (i.e., that it corresponds to a semistandard
Young tableau~$P$ of weight~$\alpha$), which is motivated only by checking one
example. We can check it for our example by applying the ``operator''~$M$
to~$\1_\alpha$ (since $\alpha=(3,4,2,3,2,2)$ is not a partition we must extend
the definitions a bit with respect to shapes of intermediate tableaux), or
equivalently by using $P_r^s=\\_{rs}$ for $r<s$ and
$P_s^s=\alpha_s-\sum_{r\in[s]}P_r^s$ for $s\in\upto6$:
$$
  \(P_r^s\)_{0\leq r\leq s<6}=
 \pmatrix{3&2&0&0&0&0\cr
           &2&1&2&0&0\cr
	   & &1&0&2&1\cr
	   & & &1&0&0\cr
	   & & & &0&1\cr
           & & & & &0\cr},
\text{whence}
  P=\Young(0,0,0,1,1|1,1,2,3,3|2,4,4,5|3|5)\,.
$$
That this condition is always satisfied is by no means clear however; in fact
it would be very difficult to associate a monomial with this property to the
sequence of operations~$e_i$ if any other criterion for selecting~$i$ (e.g.,
always taking the largest possble one) were used. If we admit that the
condition always holds, then association~II defines a map that, applied
to~$\wSem(T)$, may serve as~$\Rob_1$ (keeping in mind that Robinson works
directly with monomials rather than with semistandard Young tableaux). It is
easy to see that the shape~$\beta$ of~$M$ (i.e., the shape of the
corresponding~$P$) matches the weight of the lattice permutation found by
association~I, but not that the map~$\Rob$ defined by combining $\Rob_0$
and~$\Rob_1$ is a bijection, as required in Robinson's proof sketched above.
For this it is of vital importance that the sequence of raising operations can
be reconstructed from~$M$, which is far from obvious since the $S_{rs}$
commute, but the operations~$e_i$ do not. The key property is that in the
sequence of factors $S_{rs}$ found, the index~$s$ increases weakly, while for
fixed~$s$ the index~$r$ decreases weakly (this is the hard part); once this is
established the fact that~$M$ corresponds to a semistandard Young tableau can
also be proved easily. This is not all however, since the surjectivity
of~$\Rob$ must also be established, i.e., that every monomial~$M$ of
shape~$\beta$ arising in Young's rule applied to~$h_\alpha$, is obtained for
some word~$\wSem(T)$ of weight~$\alpha$ that corresponds under association~I
to a given lattice permutation of weight~$\beta$.

For reading words~$\wSem(P')$ of tableaux~$P'$ of \emph{partition}
shape~$\beta$, the mentioned key property is not difficult to establish;
moreover the lattice permutation is~$\wSem(\1_\beta)$ in this case, while
$P=P'$. Since we know (by theorem~\copcommthm) that the coplactic graph of any
word~$\wSem(T)$ is isomorphic to that of some such word~$\wSem(P)$ (with
$T\slid P$), we can see that the properties stated above do always hold.
However, these facts are highly non-trivial given only the information
provided in~\ref{Robinson}. Nevertheless Robinson apparently considers them as
obvious: nothing in the paper even suggests that anything needs to be proved.
In conclusion, although Robinson gives an interesting construction that
actually works, and that could be used in a proof of the Littlewood-Richardson
rule, he definitely does not provide such a proof.

\subsection Later developments.

The more recent history of the Littlewood-Richardson rule is no less
interesting than the initial phase, but since it is much more accessible and
better known, we shall limit ourselves to a brief overview. The flawed proof
given by Robinson was so incomprehensibly formulated, that its omissions
apparently went unnoticed for decades; his reasoning was reproduced
in~\ref{Littlewood} by way of proof. In the early 1960's, in~an unrelated
study, Schensted gives a combinatorial construction~\ref{Schensted} that would
later be considered to be essentially equivalent to that of Robinson; this
despite the fact that it defines a rather different kind of correspondence by
an entirely different procedure, which clearly defines a bijection, but
without obvious relation to the Littlewood-Richardson rule (although it does
involve tableaux). Nonetheless this construction would later be used in
several proofs of the rule; initially however, though the combinatorial
significance of the construction is noted by
Sch\"utzenberger~\ref{Schutzenberger QR}, this connection is not made.

This changes around 1970's, and important new properties of Schensted's
construction are found: \ref{Knuth PJM},~\ref{Greene}. From this development
emerge the first proofs of the Littlewood-Richardson rule: Lascoux and
Sch\"utzenberger introduce jeu de taquin, and using it a proof is given
in~\ref{Schutzenberger CRGS}, while Thomas gives a proof that is based
entirely on a detailed study of Schensted's construction in \ref{Thomas
thesis},~\ref{Thomas LR}. The former proof is by means of a statement similar
to our corollary~\altLRR, but it is obtained differently (e.g., for confluence
of jeu de taquin, results of \ref{Knuth PJM} and~\ref{Greene} are used). The
latter proof is interesting in that it already derives properties of
Schensted's correspondence that are related to pictures. Both approaches
differ essentially from \ref{Littlewood Richardson} and \ref{Robinson}, in
that semistandard tableaux figure in the same manner as in the present paper;
we also note that, whereas in~\ref{Robinson} the basic construction is that of
the component~$\Rob_0$ of Robinson's bijection, with many questions remaining
about~$\Rob_1$, it is $\Rob_1$ that is central in~\ref{Schutzenberger CRGS},
and $\Rob_0$ does not even occur there. Two publications from this period do
take up the construction of~\ref{Robinson}. In~\ref{Thomas Robinson} the
originally deterministic description of~$\Rob_0$ is generalised to a rewrite
system (cf.~corollary~\Robbijdescr), which is shown to be confluent.
In~\ref{Macdonald, I~(9.2)} the task of completing Robinson's argumentation
is taken up; in a long and technical proof the serious gaps in the proof
indicated above are scrupulously filled (illustrating just how much was
missing), except for the proof of the surjectivity of~$\Rob$.

After the appearance of these three proofs of the Littlewood-Richardson rule,
many publications follow; some present interesting new ideas that lead to new
proofs, but most of these are based on the same construction as one of these
earlier proofs. Of particular interest is~\ref{Zelevinsky pictures}, whose
construction relates to all three approaches. It generalises Schensted's
correspondence to pictures, which in essence consists in showing that its
bijectivity is preserved when certain restrictions parametrised by shapes
(like $\nu/\kappa$-dominance) are imposed at both sides. This is exactly
the information that is required in Schensted-based proofs (cf.~\ref{Thomas
LR}, \ref{White}, \ref{Remmel Whitney}), to make the connection with
Littlewood-Richardson tableaux. One also obtains as a special case a
bijection~$\Rob$ matching the specification~(\Robinson), which provides a
shortcut for the proof of~\ref{Macdonald} (what is not so obvious, is that
this is in fact the same correspondence as defined by~Robinson). Finally, this
establishes a symmetry between $\Rob_0$~and~$\Rob_1$ (which is what is most
notably missing in the approach of~\ref{Schutzenberger CRGS}), which allows
not only~$\Rob_1$ but also $\Rob_0$ to be defined by jeu de taquin. At the
time however, these connections were not made, and the paper does not get much
attention (this may be due to the somewhat obscure definition of pictures,
which bears no apparent relation to semistandard tableaux); these observations
are made only much later in~\ref{pictures}.

There are many more recent developments that relate to the
Littlewood-Richardson rule; these include alternative combinatorial
expressions for the coefficients~$c_{\\,\mu}^\nu$ (\ref{Berenstein
Zelevinsky}), generalisations to representations of other groups than
$\GL_n(\C)$ (\ref{Littelmann Kac-Moody}), and yet more new proofs. A
discussion of these is beyond the scope of this paper, however. We conclude by
returning to the practical computational aspects of the rule, which were the
reason it was formulated in the first place, but which seem to have moved to
the background in the course of time. To our knowledge, the first publication
to mention computer implementation of the Littlewood-Richardson rule
is~\ref{Remmel Whitney}. Surprisingly, what it calls ``a new combinatorial
rule for expanding the product of Schur functions'', is merely a translation
of the problem of multiplying $s_\\*s_\mu$ into counting objects that are
straightforward encodings of pictures $Y(\\*\mu)\to Y(\nu)$, for
varying~$\nu$; as we have indicated, this differs only marginally from the
process described in~\ref{Littlewood Richardson}. Nonetheless, the formulation
given is more straightforward to implement efficiently than most other
formulations current at that time; while the paper does not specify a detailed
algorithm, it has been used in concrete implementations (but we know of none
that are freely available at the time of writing). Currently, several freely
available and efficient implementations of the Littlewood-Richardson rule
exist in various computer algebra systems. We mention an implementation in
Maple by J.~Stembridge, contained in the package~\caps{SF}
(\URL{www.math.lsa.umich.edu/~jrs/maple.html}) and a similar implementation
in~\caps{ACE} (\URL{weyl.univ-mlv.fr/~ace/}), as well as an implementation by
the current author in the stand-alone program LiE
(\URL{wwwmathlabo.univ-poitiers.fr/~maavl/LiE/}); the latter is available for
online use, and for consultation of the documented source code~\ref{van
Leeuwen LiE}, via the mentioned \caps{WWW}-page.

\endofchapter
 \par

\Supplement References.

\input \jobname.ref

\iftoc
  \closeout\tocfile \vfil\eject
  \def\tocitem#1=#2\onpage#3.
    {\par\line{\hbox to 1.5pc{\bf #1\hss}\quad#2\dotfill#3}\ignorespaces}
  \input \jobname.toc

\fi

\bye